\documentclass[10pt,a4paper]{article}
\usepackage{fullpage,mathrsfs,graphicx,framed,color,amssymb,amsmath,amsthm}
\usepackage[boxed, noline, ruled, linesnumbered]{algorithm2e}
\usepackage{epsfig, graphicx}
\usepackage{latexsym,amsfonts,amsbsy,amssymb}
\usepackage{amsmath,amsthm}
\usepackage{enumerate}
\usepackage{bm}
\makeatletter 

\@addtoreset{equation}{section}
\makeatother 
\allowdisplaybreaks 
\newtheorem{theorem}{Theorem}[section]
\newtheorem{lemma}{Lemma}[section]
\newtheorem{cor}{Corollary}[section]
\newtheorem{remark}{Remark}[section]

\begin{document}
\title{Sequential symmetric interior penalty discontinuous Galerkin method for fully coupled quasi-static thermo-poroelasticity problems}
\author{Fan Chen\footnote{Department of Mathematics, Beijing University of Technology, Beijing 100124, China.
({\tt ChenFan2021@emails.bjut.edu.cn})},\ \ \
Ming Cui\footnote{Department of Mathematics, Beijing University of Technology, Beijing 100124, China. ({\tt mingcui@bjut.edu.cn})} \  \  and\  \ Chenguang Zhou\footnote{Department of Mathematics, Beijing University of Technology, Beijing 100124, China. (Corresponding author: {\tt zhoucg@bjut.edu.cn})} }
\date{}
\maketitle
\begin{abstract}
In this paper, we investigate a sequentially decoupled numerical method for solving the fully coupled quasi-static thermo-poroelasticity problems with nonlinear convective transport. The symmetric interior penalty discontinuous Galerkin method is employed for spatial discretization and the backward Euler method for temporal discretization. Unlike other splitting algorithms, this type of sequential method does not require any internal iterations and the computational efficiency is higher than that of the fully implicit nonlinear numerical scheme. In the theoretical analysis, a cut-off operator is introduced to prove the existence and uniqueness of numerical solution and the stability analysis of numerical scheme is conducted. Then, we derive the optimal convergence order estimates in space and time. Finally, several numerical examples are presented to illustrate the accuracy and efficiency of our proposed method.

\end{abstract}

\section{Introduction}\label{intro}
Poroelasticity problems, which are also commonly referred to as the Biot model, describe the interaction between solid deformation and
fluid flow in porous media, and the theories go back to the works of Terzhagi \cite{Terzaghi1943Theoretical} and Biot \cite{12}. The Biot model has been treated by various numerical methods, such as finite difference methods \cite{MR1930843,MR1957690}, finite volume methods \cite{MR2280930,MR3478962}, finite element methods \cite{MR4684196,MR2177147,MR4931528}, virtual element methods \cite{MR4636155,MR4647655,MR4486103}, discontinuous Galerkin (DG) methods \cite{MR4734247,33,19} and weak Galerkin methods \cite{MR4588484,MR4799262,MR4136610}. The thermo-poroelasticity model is an extension of Biot model to the non-isothermal case. Due to the generality and strong physical grounding, the thermo-poroelasticity framework is widely used in various fields, such as simulating soil consolidation processes and energy transfer in geotherm \cite{poromechanics}, thermal ablation of tumors in biomechanics \cite{nzl}, and carbon dioxide sequestration and thermal remediation in environmental science \cite{2017Multi}.

The thermo-poroelasticity model is firstly established in \cite{MR3825656} by Brun et al., using a two-scale expansion method to describe the control equation system of fully saturated fine-grained porous media under quasi-static deformation conditions. The system includes a pore fluid flow equation, an elasticity equation and a heat conduction equation, representing mass conservation, momentum conservation and energy conservation, respectively. According to the constitutive laws and nonlinearities considered, the system can be linear or nonlinear. In this work, we mainly consider the nonlinear thermo-poroelasticity model and design stable and efficient numerical methods.

Regarding the research of theoretical analysis and numerical methods of nonlinear thermo-poroelasticity model, Brun et al. \cite{2018Well} consider the well-posedness of nonlinear thermo-poroelasticity problem and provide the a priori energy estimates and the regularity properties of solutions. In \cite{2019Monolithic}, Brun et al. present the monolithic and splitting-based iterative procedures for the coupled nonlinear thermo-poroelasticity problem, which transforms a nonlinear and fully coupled problem into a set of simpler subproblems to be solved sequentially in an iterative fashion. In \cite{MR4432106}, Zhang and Rui study the standard Galerkin method for fully coupled nonlinear thermo-poroelasticity problems and analyze the error convergence order. In \cite{MR4579739}, Antonietti et al. concentrate on the four-field nonlinear thermo-poroelasticity problem by taking into account one extra variable and a priori $hp$-version error estimates of the linearized DG semi-discrete scheme are derived. In \cite{chen}, Chen et al. provide a fully implicit nonlinear discrete scheme, which is constructed by adopting the DG method for the spatial approximation and the backward Euler method for the temporal discretization.

Due to the complexity of the coupled multi-physical field thermo-poroelastic model, it is very important to develop an accurate and computationally efficient numerical method. To solve the thermo-poroelastic model, there are mainly three strategies: full implicit method, iterative decoupling method and sequential decoupling method. The fully implicit method \cite{MR3899054,MR1690489} couples and solves all variables at each time step, requiring the construction and solution of a high-dimensional algebraic system. Although fully implicit methods have good stability, their high computational cost limits the application in large-scale engineering problems. The iterative decoupling method \cite{MR4685928,MR3614285,2019Monolithic,MR3253752} iteratively solves the subproblem based on the approximate solution of the previous iteration at each time step until the preset convergence tolerance is met. The method can effectively reduce the computation complexity of single step, but the iterative process may incur additional computational costs, especially in strongly nonlinear situations where the convergence is difficult to guarantee. The sequential decoupling method \cite{MR3800024,MR4101493} successively solves sub-models at each time step and does not require iteration. Compared with the fully implicit method, the dimension of discrete system is reduced and the computational efficiency is improved. Recently, Chaabane and Rivi\`ere \cite{MR3742883} propose a novel sequential discontinuous Galerkin (SDG) method for solving Biot systems, which separates the coupled models and solves two subproblems sequentially at each time step. By adding the stabilization terms to ensure the stability of numerical scheme, the a priori error estimates with respect to displacement and pressure are obtained. Subsequently, Chaabane and Rivi\`ere extend this type of sequential idea to the case of conforming finite element method \cite{MR3777103}, and in \cite{MR4457768}, Shen et al. make use of the SDG method for solving two-phase flow poroelasticity equations. However, the above-mentioned three papers only prove the convergence order estimates that are optimal in space and suboptimal in time. Furthermore, no stability analysis is conducted.

In this paper, we propose a sequential decoupling method for the fully coupled quasi-static thermo-poroelasticity problems with nonlinear convective transport. This method is divided into three steps at each time step: In the first step, assuming that the temperature and displacement of the first two time layers are given, the fluid equation is used to solve the pressure; In the second step, the values of the first two time layers of known displacement are used to solve the temperature by utilizing the obtained pressure and the heat equation; In the third step, the displacement is solved by applying the obtained pressure, temperature and the mechanical equation. In this way, we have fully decoupled the thermo-poroelasticity model at each time step. It can be seen that the above method process relies on the pressure, temperature and displacement of the first two time layers. Here, we provide two ways to obtain them. The first way follows the idea proposed in \cite{MR3742883} and does not require rewriting the program code. Since we remove the divergence and stabilization terms in this strategy, it is necessary to assume that the initial time step size is very small in order to ensure the change of displacement with respect to time very small, and even can be ignored. The second way selects a fully implicit nonlinear numerical scheme in \cite{chen} to obtain the optimal convergence without the assumption that the initial time step size is very small. In terms of space and time discretization, the DG method and backward Euler method are adopted. We analyze the existence and uniqueness of numerical solution, and the stability of algorithm. Then, the optimal convergence order estimates in space and time are derived. To the best of our knowledge, this is the first paper to investigate the stability of SDG method presented in \cite{MR3742883} and obtain the optimal error estimates in both space and time.

The outline of this paper goes as follows. In Section \ref{EQ}, we introduce the quasi-static thermo-poroelasticity model and provide several necessary notations. In Section \ref{DG}, we give the details of sequential decoupling algorithm, and conduct the well-posedness analysis of algorithm. After the definition of some projection operators, the optimal convergence order estimates are derived in Section \ref{EA}. And finally, we report numerical experiments in Section \ref{NE} to validate the accuracy and efficiency of our proposed method.

\section{Mathematical model}\label{EQ}
The following fully coupled quasi-static thermo-poroelasticity problem with nonlinear convective transport is the main focus of this work. Let $\Omega$  be a convex polygonal or polyhedral domain in $\mathbb R^d$ $(d=2,3)$ with Lipschitz boundary $\partial\Omega$, and $(0,t_f]$ a time interval with $t_f>0$. Given a body force $\bm f$, a mass source $g$ and a heat source $z$, seek the displacement ${\mathbf u}(t):\Omega\to\mathbb R^d$, the pressure $p(t):\Omega\to\mathbb R$ and the temperature $T(t):\Omega\to\mathbb R$ such that
\begin{align}
 -\nabla\cdot\bm\sigma({\mathbf u})+\alpha\nabla p+\beta\nabla T
&=\bm{f}, \ \ \textrm {in} \  \Omega\times(0,t_f], \label{e1}\\
\frac {\partial }{\partial t}(c_0p-b_0T+\alpha\nabla\cdot{\mathbf u})-\nabla\cdot({\mathbf K}\nabla p)&=g, \ \  \textrm {in}  \ \Omega\times(0,t_f],\label{e2}\\
\frac {\partial }{\partial t}(a_0T-b_0p+\beta\nabla\cdot{\mathbf u})-{\mathbf K}\nabla p\cdot\nabla T-\nabla\cdot(\bm\Theta\nabla T)&=z, \ \  \textrm {in}  \ \Omega\times(0,t_f],\label{e3}
\end{align}
with the boundary conditions
\begin{align}
{\mathbf u}=\mathbf 0,\ \ p=0,\ \ T&=0, \ \  \textrm {on} \ \partial\Omega\times(0,t_f],\label{e4}
\end{align}
and the initial conditions
\begin{align}
{\mathbf u}(\cdot,0)&={\mathbf u}^0, \ \ \textrm {in} \ \Omega,\label{e5}\\
p(\cdot,0)&=p^0, \ \ \textrm {in} \ \Omega,\label{e6}\\
T(\cdot,0)&=T^0, \ \ \textrm {in} \ \Omega.\label{e7}
\end{align}
Here, the effective stress tensor is denoted by $\bm\sigma({\mathbf u})=2\mu\bm\epsilon({\mathbf u})+\lambda\nabla\cdot{\mathbf u\mathbf I}$, where  $\bm\epsilon({\mathbf u})=\frac{\nabla{\mathbf u}+\nabla{\mathbf u}^\mathsf{T}}{2}$ is the strain tensor, ${\mathbf I}$ is the identity tensor, $\mu$ and $\lambda$ are the $Lam\acute{e}$ constants. $\alpha$ and $\beta$ are the Biot-Willis constant and the thermal stress coefficient, respectively. $c_0 \geq 0$ represents the constrained specific storage coefficient. The effective volumetric heat capacity is $a_0$, and the thermal dilation coefficient is $b_0$. ${\mathbf K}=(K_{ij})_{i,j=1}^d$ is the permeability divided by fluid viscosity, and $\bm\Theta=(\Theta_{ij})_{i,j=1}^d$ is  the effective thermal conductivity.

Now, define some function spaces that will be used throughout this article. We use the conventional notion of Sobolev space $W^{k,l}(\Omega)$ (cf. \cite{MR0450957}) in $L^l(\Omega)$, which admits weak derivatives up to order $k$ in the same space. Specifically, we indicate for $k\geq 0$ by $H^{k}(\Omega):=W^{k,2}(\Omega)$. The associated inner-product and norm (resp.seminorm) in $H^k(\Omega)$ are denoted by $(\cdot, \cdot)_k$ and $\|\cdot\|_k$ (resp. $|\cdot|_k$), respectively. When $k=0$, $H^0(\Omega)$ coincides with the space of square-integrable functions $L^2(\Omega)$. In this case, the subscript $k$ is suppressed from the notation of inner product and norm. In addition, on subregion $D\subset\Omega$, the norm of $L^2(D)$ is denoted by $\|\cdot\|_{0,D}$. Denote by $H_0^1(\Omega)$ the subspace of $ H^1(\Omega)$ consisting of the functions with vanishing trace on $\partial\Omega$. Denote by $||\cdot||_{0,\infty}$ the norm on $L^{\infty}(\Omega)$.

To deal with the functions involving time and space, we consider the space
$$ L^2(0,t_f;H^k(\Omega))=\{u:(0,t_f]\to H^k(\Omega),\int_0^{t_f}||u||_k^2dt<\infty\},$$
with the norm
$$||u||_{L^2(0,t_f;H^k(\Omega))}^2=\int_0^{t_f}||u||_{k}^2dt.$$

Before providing the numerical scheme of the problem  \eqref{e1}-\eqref{e7}, the following  assumptions are based on \cite{2018Well}.

{\bf Assumption  A.} ({\bf A1}) The permeability ${\mathbf K}$ and the heat conductivity $\bm\Theta$ are symmetric and positive definite, and ${\mathbf K}$, $\bm\Theta\in(L^{\infty}(\Omega))^{d\times d}$. Additionally, there exist positive constants $k_m$ and $k_M$ such that
\begin{equation*}
k_m|\zeta|^2\leq \zeta^\mathsf{T}{\mathbf K}(\mathbf x)\zeta, \ \ |{\mathbf K}(\mathbf x)\zeta| \leq k_M|\zeta|,\ \ \forall \zeta\in \mathbb R^d,\ \ \mathbf x \in \Omega,
\end{equation*}
and there exist constants $\theta_m>0$ and $\theta_M>0$ such that
\begin{equation*}
\theta_m|\zeta|^2\leq \zeta^\mathsf{T}\bm\Theta(\mathbf x)\zeta,\ \ |\bm\Theta(\mathbf x)\zeta|\leq \theta_M|\zeta|,\ \ \forall \zeta\in \mathbb R^d,\ \ \mathbf x \in \Omega.
\end{equation*}
({\bf A2}) The parameters $a_0$, $b_0$, $c_0$, $\alpha$, $\beta$, $\mu$ and $\lambda$ are strictly postive constants.\\
({\bf A3}) The parameters $a_0$, $b_0$ and $c_0$ satisfy $c_0-2b_0>0,a_0-2b_0\ge0$.\\
({\bf A4}) The source terms satisfy  $z,g\in L^2(0,t_f;L^2(\Omega))$ and ${\bm f}\in H^1(0,t_f;(L^2(\Omega))^d)$.\\
({\bf A5}) The initial conditions satisfy $p^0,T^0\in H_0^1(\Omega)$ and ${\mathbf u}^0\in(H_0^1(\Omega))^d$.\\
({\bf A6}) For $l>d$ and $k_i\geq1$ $(i=1,2,3)$, the functions $\nabla p$, $p$, $T$, $\mathbf u$ satisfy
\begin{align*}
&||\nabla p||_{L^{\infty}(0,t_f;W^{1,l}(\Omega))}+||p||_{L^{\infty}(0,t_f;W^{2,l}(\Omega))}\\
+&||p||_{L^{\infty}(0,t_f;H^{k_2+1}(\Omega))}+||T||_{L^{\infty}(0,t_f;H^{k_3+1}(\Omega))}+
||\mathbf u||_{L^{\infty}(0,t_f;(H^{k_1+1}(\Omega))^d)}\\
+&||p_t||_{L^{\infty}(0,t_f;W^{1,l}(\Omega))}+||T_t||_{L^{\infty}(0,t_f;W^{1,l}(\Omega))}+
||\mathbf u_t||_{L^{\infty}(0,t_f;(W^{2,l}(\Omega))^d)}\\
+&||p_{tt}||_{L^{2}(0,t_f;L^{l}(\Omega))}+||T_{tt}||_{L^{2}(0,t_f;L^{l}(\Omega))}+
||\mathbf u_{tt}||_{L^{2}(0,t_f;(W^{1,l}(\Omega))^d)}\le C.
\end{align*}

Furthermore, the symbol $C$, either with or without subscripts, denotes a generic positive constant which may be different at its different occurrences throughout the paper.

\section{Sequential decoupling algorithm}\label{DG}
\subsection{Algorithm details}\label{VF}
Let $\mathcal T_h$ denote a shape-regular triangulation of the domain $\Omega$, and $K$ represents a triangle element with the diameter $h_K$ and $h=\max\limits_{K\in \mathcal T_h} h_K$. We define all element edges by $\Gamma$ and $h_e$ stands for the length of an edge $e\in\Gamma$. On account of the polygonal or polyhedral domain $\Omega$, we have $\Gamma=\Gamma_I \cup \partial\Omega,$ where $\Gamma_I$ includes all the internal edges.

We define the following broken Sobolev space with integer $k\geq 0$,
$$H^k(\mathcal{T}_h)=\{v\in L^2(\Omega):v|_K\in H^k(K),\ \forall K\in\mathcal{T}_h\}.$$

Now, we introduce the average and jump operators that are required for the DG method. Let $K_i$ and $K_j$ $(i>j)$ be two adjacent elements of $\mathcal T_h$ which share a common edge $e$, i.e., an interior edge $e=\partial K_i\cap\partial K_j\subset\Gamma_I.$  We assume that the unit normal vector ${\mathbf n_e}$ is oriented from $K_i$ to $K_j$. The jump and average of $v$ on $e$ are given by
\[
[v]=v|_{K_i}-v|_{K_j},\   \   \{v\}=\frac{1}{2}(v|_{K_i}+v|_{K_j}).
\]
If we have a domain boundary edge $e\subset\partial K_i\cap\partial\Omega$, then
$$[v]=v|_{e},\   \   \{v\}=v|_{e}.$$

Let $k_i\geq1$ $(i=1,2,3)$ be a positive integer, we define the spaces of discontinuous finite element as follows,
\begin{align*}
{\mathbf V}_h&=\{{\mathbf v}\in (L^{2}(\Omega))^d:{\mathbf v}|_{K}\in (P_{k_1}(K))^d,\ \forall K\in \mathcal T_h\},\\
Q_h&=\{q\in L^{2}(\Omega):q|_{K}\in P_{k_2}(K),\ \forall K\in\mathcal T_h\},\\
S_h&=\{w\in L^{2}(\Omega):w|_{K}\in P_{k_3}(K),\ \forall K\in\mathcal T_h\},
 \end{align*}
where $P_{k_i}(K)$ is the space spanned by the polynomials of degree at most $k_i$ on each element $K$.

For the convenience of illustrating algorithm, we introduce several notations, for any ${\mathbf v}$, ${\mathbf w} \in {\mathbf V}_h$ and any $p$, $q\in Q_h$ (or $S_h$),
$$\begin{array}{lll}
a({\mathbf v},{\mathbf w})=&\sum\limits_{K\in \mathcal T_h}\int_{K}\bm\sigma({\mathbf v}):\bm\epsilon({\mathbf w})dK-\sum\limits_{e\in \Gamma}\int_{e}\{\bm\sigma({\mathbf v})\}{\bm n_e}\cdot[{\mathbf w}]ds-\sum\limits_{e\in \Gamma}\int_{e}\{\sigma({\mathbf w})\}{\bm n_e}\cdot[{\mathbf v}]ds\\
&+\sigma_1\sum\limits_{e\in \Gamma}\int_{e}h_e^{-1}[{\mathbf v}]\cdot[{\mathbf w}]ds,
\end{array}
$$
$$b({\mathbf v},p)=-\sum\limits_{K\in \mathcal T_h}\int_K \nabla p\cdot{\mathbf v}dK+\sum\limits_{e\in \Gamma_I}\int_{e}\{{\mathbf v}\}\cdot{\bm n_e}[p]ds=\sum\limits_{K\in \mathcal T_h}\int_K p\nabla\cdot{\mathbf v}dK-\sum\limits_{e\in \Gamma_I}\int_{e}\{p\}[{\mathbf v}]\cdot{\bm n_e}ds,$$
and
$$c(\phi;p,q)=\sum\limits_{K\in \mathcal T_h}\int_K\phi\nabla p\cdot\nabla qdK-\sum\limits_{e\in \Gamma}\int_{e}\{\phi\nabla p\}\cdot{\bm n_e}[q]ds-\sum\limits_{e\in \Gamma}\int_{e}\{\phi\nabla q\}\cdot{\bm n_e}[p]ds+\sigma_2\sum\limits_{e\in \Gamma}\int_{e}h_e^{-1}[p][q]ds,$$
where $\sigma_1$ and $\sigma_2$ are two large enough positive constants, and $\phi$ is a symmetric positive definite matrix. Moreover, there exist constants $\phi_m>0$ and $\phi_M>0$ such that
\begin{equation*}
\phi_m|\zeta|^2\leq \zeta^\mathsf{T}\phi(\mathbf x)\zeta, \ \ |\phi(\mathbf x)\zeta| \leq \phi_M|\zeta|,\quad \forall \zeta\in \mathbb{R}^d,\ \ \mathbf x\in \Omega.
\end{equation*}

For temporal discretization, let $0 = t_0 < t_1 < t_2 < \cdots < t_N = t_f$ be a partition of the time interval $(0,t_f]$. We choose two time step sizes $\tau_0$ and $\tau$ $(\tau_0\ll\tau)$ for reasons that will be apparent in Remark \ref{remark1}, and define
\begin{align*}
t_1&=\tau_0,\\
t_n&=t_1+(n-1)\tau,\ \ \textrm {for any integer } n\geq2.
\end{align*}
We indicate the approximation of ${\mathbf u}(t_n)$, $p(t_n)$ and $T(t_n)$ by the notations ${\mathbf u}_h^n$, $p_h^n$ and $T_h^n$, respectively. By ${\bm f}^n$, $g^n$ and $z^n$, we denote ${\bm f}(t_n)$, $g(t_n)$ and $z(t_n)$, separately.

Before giving the SDG scheme, we introduce the cut-off
operator $\mathcal M$ as described in \cite{2019Monolithic,MR2039576,MR2116915},
\begin{equation}\label{m}
\mathcal M({\mathbf z})(x):=\left\{ \begin{array}{l}
{\mathbf z}(x) , \hspace{2cm} |{\mathbf z}(x)|\le M,\\
M{\mathbf z}(x)/|{\mathbf z}(x)|, \hspace{0.6cm} |{\mathbf z}(x)|> M,
\end{array} \right.
\end{equation}
where $M$ is a large positive constant. It should be noted that this operator can help to facilitate the convergence analysis in the following algorithm. Obviously, if the exact Darcy flux is bounded, i.e., ${\mathbf K}\nabla p^n \in (L^{\infty}(\Omega))^d$, and select $M$ large enough, we have $\mathcal M({\mathbf K}\nabla p^n)(x)={\mathbf K}\nabla p^n(x)$. Thus, a precise value for the constant $M$ is not necessary.

Now, we are ready to define the SDG method of \eqref{e1}-\eqref{e7}, which consists of three steps at each time step for $n\geq1$,

{\bf Step1.} Given $p_h^n\in Q_h$, $T_h^n$, $T_h^{n-1}\in S_h$ and
${\mathbf u}_h^n$, ${\mathbf u}_h^{n-1}\in {\mathbf V}_h$, find $p_h^{n+1}\in Q_h$ such that
\begin{align}\label{s1}
&c_0(\frac{p_h^{n+1}-p_h^{n}}{\tau},q_h)-b_0(\frac{T_h^{n}-T_h^{n-1}}{\tau},q_h)+\alpha b(\frac{{\mathbf u}_h^{n}-{\mathbf u}_h^{n-1}}{\tau},q_h)\notag\\
&+c({\mathbf K};p_h^{n+1},q_h)=(g^{n+1},q_h),\ \  \forall q_h\in Q_h.
\end{align}

{\bf Step2.} Given $p_h^n$, $p_h^{n+1}\in Q_h$, $T_h^n\in S_h$ and
${\mathbf u}_h^n$, ${\mathbf u}_h^{n-1}\in {\mathbf V}_h$, find $T_h^{n+1}\in S_h$ such that
\begin{align}\label{s2}
&a_0(\frac{T_h^{n+1}-T_h^{n}}{\tau},s_h)-b_0(\frac{p_h^{n+1}-p_h^{n}}{\tau},s_h)+\beta b(\frac{{\mathbf u}_h^{n}-{\mathbf u}_h^{n-1}}{\tau},s_h)\notag\\
&-(\mathcal M(\mathbf K\nabla p_h^{n})\cdot\nabla T_h^{n+1},s_h)+c(\bm\Theta;T_h^{n+1},s_h)=(z^{n+1},s_h),\ \ \forall s_h\in S_h.
\end{align}

{\bf Step3.} Given $p_h^{n+1}\in Q_h$, $T_h^{n+1}\in S_h$ and
${\mathbf u}_h^n$, ${\mathbf u}_h^{n-1}\in {\mathbf V}_h$, find ${\mathbf u}_h^{n+1}\in {\mathbf V}_h$ such that
\begin{align}\label{s3}
&a({\mathbf u}_h^{n+1},{\mathbf v_h} )-\alpha b({\mathbf v}_h,p_h^{n+1})-\beta b({\mathbf v}_h,T_h^{n+1})\notag\\
&+\gamma(\frac{{\mathbf u}_h^{n+1}-{\mathbf u}_h^{n}}{\tau},{\mathbf v_h})-\gamma(\frac{{\mathbf u}_h^{n}-{\mathbf u}_h^{n-1}}{\tau},{\mathbf v_h})=({\bm f}^{n+1},{\mathbf v_h}),\ \ \forall {\mathbf v}_h\in {\mathbf V}_h,
\end{align}
where $(\cdot,\cdot)=\sum\limits_{K\in\mathcal T_h}(\cdot, \cdot)_K$, and $\gamma$ is a stabilization parameter \cite{MR3742883}. In the above method, we utilize $(\mathcal M(\mathbf K\nabla p_h^{n})\cdot\nabla T_h^{n+1},s_h)$ to approximate the convective coupling term instead of $(\mathbf K\nabla p_h^{n}\cdot\nabla T_h^{n+1},s_h).$ It is shown in \cite{MR3904522} for the poroelasticity model that if the flux is bounded on the continuous level, the discrete flux will inherit this property. Then, with sufficient regular domain, source terms and initial data, the existence of the constant $M$ in \eqref{m} is guaranteed \cite{2019Monolithic}.

In order to start the algorithm, we need to address the issue of approximating $p_h^1$, $T_h^1$ and ${\mathbf u}_h^1$. The initial values are chosen to be the elliptic projections of the initial data,
\begin{equation*}
{\mathbf u}_h^0=R_u{\mathbf u}^0,\ \ p_h^0=R_pp^0,\ \ T_h^0=R_TT^0,
\end{equation*}
where the projection operators $R_u$, $R_p$ and $R_T$ are defined as follows: Define the projection $R_u{\mathbf u}\in {\mathbf V}_h$ of displacement satisfying
\begin{equation}
a({\mathbf u}-R_u{\mathbf u}, {\mathbf v}_h)=0,\ \ \forall {\mathbf v}_h\in {\mathbf V}_h.\label{ea1}
\end{equation}
Define the projection $R_hp\in Q_h$ of pressure satisfying
\begin{equation}
c(\mathbf{K};p-R_hp, q_h)=0,\ \ \forall q_h\in Q_h.\label{ea2}
\end{equation}
Define the projection $R_hT\in S_h$ of temperature satisfying
\begin{equation}
c(\bm\Theta;T-R_hT, s_h)-({\mathbf K}\nabla p\cdot\nabla(T-R_hT),s_h)=0,\ \ \forall s_h\in S_h.\label{ea3}
\end{equation}
And $(p_h^1,T_h^1,{\mathbf u}_h^1)$ can be solved by applying the following two options.

\textbf{Option 1}: For the convenience of writing program code, we follow the idea in \cite{MR3742883} and calculate $(p_h^1, T_h^1, {\mathbf u}_h^1)$ in the following three steps:

{\bf Initial Step1.} Find $p_h^{1}\in Q_h$ such that
\begin{equation}\label{is1}
c_0(\frac{p_h^{1}-p_h^{0}}{\tau_0},q_h)+c({\mathbf K};p_h^{1},q_h)=(g^{1},q_h),\ \ \forall q_h\in Q_h.
\end{equation}

{\bf Initial Step2.} Find $T_h^{1}\in S_h$ such that
\begin{equation}\label{is2}
a_0(\frac{T_h^{1}-T_h^{0}}{\tau_0},s_h)-b_0(\frac{p_h^{1}-p_h^{0}}{\tau_0},s_h)-(\mathcal M(\mathbf K \nabla p_h^{0})\cdot\nabla T_h^{1},s_h)+c(\bm\Theta;T_h^{1},s_h)=(z^{1},s_h),\ \ \forall s_h\in S_h.
\end{equation}

{\bf Initial Step3.} Find ${\mathbf u}_h^{1}\in {\mathbf V}_h$ such that
 \begin{equation}\label{is3}
a({\mathbf u}_h^{1},{\mathbf v_h} )-\alpha b({\mathbf v}_h,p_h^{1})-\beta b({\mathbf v}_h,T_h^{1})=({\bm f}^{1},{\mathbf v_h}),\ \ \forall{\mathbf v}_h\in {\mathbf V}_h.
\end{equation}

\begin{remark}\label{remark1}
Compared with the scheme \eqref{s1}-\eqref{s3}, the initial scheme \eqref{is1}-\eqref{is3} leaves out the term $b_0(\frac{T_h^{n}-T_h^{n-1}}{\tau},q_h),$ the divergence terms $\alpha b(\frac{{\mathbf u}_h^{n}-{\mathbf u}_h^{n-1}}{\tau},q_h)$, $\beta b(\frac{{\mathbf u}_h^{n}-{\mathbf u}_h^{n-1}}{\tau},s_h)$ and the stabilization terms $\gamma(\frac{{\mathbf u}_h^{n+1}-{\mathbf u}_h^{n}}{\tau},{\mathbf v_h})$, $\gamma(\frac{{\mathbf u}_h^{n}-{\mathbf u}_h^{n-1}}{\tau},{\mathbf v_h})$. Here, we assume that by selecting a small initial time step size $\tau_0$, the change in time of $\nabla \cdot \mathbf u$ should be small enough to be neglected.
\end{remark}

\textbf{Option 2}: We adopt a fully implicit numerical scheme to obtain $(p_h^1,T_h^1,{\mathbf u}_h^1)$, which is the optimal convergence \cite{chen}. The numerical scheme reads: Seek $({\mathbf u}_h^1, p_h^1,T_h^1)\in {\mathbf V}_h\times Q_h\times S_h$  such that
\begin{align}
&a({\mathbf u}_h^1,{\mathbf v_h} )-\alpha b({\mathbf v}_h,p_h^1)-\beta b({\mathbf v}_h,T_h^1)=({\bm f}^1,{\mathbf v_h}),\ \ \forall {\mathbf v}_h\in {\mathbf V}_h,\label{b71}\\
&c_0(\frac{p_h^{1}-p_h^{0}}{\tau_0},q_h)-b_0(\frac{T_h^{1}-T_h^{0}}{\tau_0},q_h)+\alpha b(\frac{{\mathbf u}_h^{1}-{\mathbf u}_h^{0}}{\tau_0},q_h)+c({\mathbf K};p_h^1,q_h)=(g^1,q_h),\ \ \forall q_h\in Q_h,\label{b81}\\
&a_0(\frac{T_h^{1}-T_h^{0}}{\tau_0},s_h)-b_0(\frac{p_h^{1}-p_h^{0}}{\tau_0},s_h)+\beta b(\frac{{\mathbf u}_h^{1}-{\mathbf u}_h^{0}}{\tau_0},s_h)+c(\bm\Theta;T_h^1,s_h)\notag\\
&\hspace{7.1cm}-({\mathbf K}\nabla p_h^{1}\cdot\nabla T_h^1,s_h)=(z^1,s_h),\ \ \forall s_h\in S_h.\label{b91}
\end{align}
Here the time step size $\tau_0$ does not need to be very small, and we can choose $\tau_0=\tau$.

\subsection{Well-posedness}\label{EU}
The following semi-norms and norms are introduced for theoretical analysis:
\begin{align*}
&|\!|\!|{\mathbf v}|\!|\!|_V^2=\sum\limits_{K\in \mathcal T_h}||{\bm\epsilon(\mathbf v)}||_{0,K}^2+\sum\limits_{e\in \Gamma}\int_e h_e^{-1}[{\mathbf v}]^2ds=||{\bm\epsilon(\mathbf v)}||^2
+|{\mathbf v}|_\ast^2,\ \ \forall{\mathbf v}\in {\mathbf V}_h,\\
&|\!|\!|q|\!|\!|^2=\sum\limits_{K\in \mathcal T_h}||\nabla q||_{0,K}^2+\sum\limits_{e\in \Gamma}\int_eh_e^{-1}[q]^2ds=|q|_1^2+|q|_\ast^2,\ \ \forall q\in Q_h\ (or\ S_h).
\end{align*}

Before the well-posedness analysis of SDG method, we need to supply the following
useful lemmas.
\begin{lemma}\label{lem2.2}\cite{25} Let $e$ denote an edge of the element $K\in \mathcal T_h$. Then, there exist positive constants $C$ independent of the mesh size $h$ such that
\begin{align}
||v||_{0,e}&\leq Ch_K^{-\frac{1}{2}}||v||_{0,K},\ \ \forall v\in P_k(K),\\
||\nabla v\cdot\bm n_e||_{0,e}&\leq Ch_K^{-\frac{1}{2}}||\nabla v||_{0,K},\ \ \forall v\in P_k(K),\\
||v||_{0,e}&\leq Ch_K^{-\frac{1}{2}}(||v||_{0,K}+h_K||\nabla v||_{0,K}),\ \ \forall v\in H^1(K),\\
||\nabla v\cdot\bm n_e||_{0,e}&\leq Ch_K^{-\frac{1}{2}}(||\nabla v||_{0,K}+h_K||\nabla^2 v||_{0,K}),\ \ \forall v\in H^2(K).
\end{align}
\end{lemma}

\begin{lemma}\label{lem2.21}\cite{25} Let $e$ denote an edge of the element $K\in \mathcal T_h$. Then we have
\begin{align}
||\nabla^jv||_{0,K}&\leq Ch_K^{-j}||v||_{0,K},\ \  0\le j\le k.
\end{align}
\end{lemma}

From \cite{MR2519594}, it follows that the Korn's inequality below.
\begin{lemma}\label{lem2.3}\cite{6} There exists a positive constant $C$ independent of the mesh size $h$, such that
\begin{equation*}
\sum\limits_{K\in \mathcal T_h}||\nabla\mathbf v||_{0,K}^2+\sum\limits_{e\in \Gamma}\int_e h_e^{-1}[{\mathbf v}]^2ds\leq C|\!|\!|{\mathbf v}|\!|\!|_V^2,\ \ \forall {\mathbf v}\in {\mathbf V}_h.
\end{equation*}
\end{lemma}

\begin{lemma}\label{lem2.4}\cite{6} There exists a constant $C>0$ such that
\begin{align*}
&||{\mathbf v}||^2\leq C|\!|\!|{\mathbf v}|\!|\!|_V^2,\ \ \forall {\mathbf v}\in {\mathbf V}_h,\\
&||q||^2\leq C|\!|\!|q|\!|\!|^2,\ \ \forall q\in Q_h \ (or\ S_h).
\end{align*}
\end{lemma}

\begin{lemma}\label{lem2.5}\cite{25} If $\sigma_1$ is large enough, there exist positive constants $C_1$ and $C_2$, independent of the mesh size $h$, such that
\begin{align*}
&a({\mathbf v},{\mathbf v})\geq C_1|\!|\!|{\mathbf v}|\!|\!|_V^2,\ \ \forall{\mathbf v}\in{\mathbf V}_h,\\
&a({\mathbf v},{\mathbf w})\leq C_2|\!|\!|{\mathbf v}|\!|\!|_V|\!|\!|{\mathbf w}|\!|\!|_V,\ \ \forall{\mathbf v},{\mathbf w}\in{\mathbf V}_h.
\end{align*}
Similarly, if $\sigma_2$ is large enough, there exist positive constants $C_3$ and $C_4$, independent of $h$, such that
\begin{align*}
&c(\phi;p,p)\geq C_3|\!|\!|p|\!|\!|^2,\ \ \forall p\in Q_h \ (or\ S_h),\\
&c(\phi;p,q)\leq C_4|\!|\!|p|\!|\!||\!|\!|q|\!|\!|,\ \ \forall p,q\in Q_h\ (or\ S_h),
\end{align*}
where the matrix $\phi$ is symmetric and positive definite.
\end{lemma}

Next, we prove the existence and uniqueness of the SDG method. To this aim, we show that the cut-off operator $\mathcal M$ is uniformly Lipschitz continuous in the following sense.
\begin{lemma}\label{lem0}\cite{MR2116915} (Property of $\mathcal M$). The cut-off operator $\mathcal M$ defined in \eqref{m} is uniformly Lipschitz continuous, i.e.,
\begin{equation*}
||\mathcal M(\mathbf z_1)-\mathcal M(\mathbf z_2)||_{0,\infty}\le||\mathbf z_1-\mathbf z_2||_{0,\infty}.
\end{equation*}
Then, we have
\begin{equation*}
||\mathcal M({\mathbf K}\nabla p^n)-\mathcal M({\mathbf K}\nabla p_h^n)||_{0,\infty}\le||{\mathbf K}\nabla p^n-{\mathbf K}\nabla p_h^n||_{0,\infty},
\end{equation*}
and
\begin{equation*}
||\mathcal M({\mathbf K}\nabla p_h^n)||_{0,\infty}\le M.
\end{equation*}
\end{lemma}

\begin{theorem}\label{eu}
Assuming that (A1)-(A6) hold true, then when $\tau \le\frac{2C_3a_0}{M^2}$, the solution $({\mathbf u}_h^n,p_h^n,T_h^n)\in {\mathbf V}_h\times Q_h\times S_h$ of \eqref{s1}-\eqref{s3} exists and is unique.
\end{theorem}
\begin{proof}
\indent Since \eqref{s1}-\eqref{s3} is linear and finite dimensional, the existence is equivalent to the uniqueness. For $n\geq 2$ and assuming that $\widetilde{p}_h$ is the difference between two solutions of \eqref{s1}, we have
\begin{equation*}
c_0(\frac{\widetilde{p}_h}{\tau},q_h)+c({\mathbf K};\widetilde{p}_h,q_h)=0.
\end{equation*}
Choosing $q_h=\widetilde{p}_h$ in the equation above and applying Lemma \ref{lem2.5}, we have
\begin{equation*}
\frac{c_0}{\tau}||\widetilde{p}_h||^2+C_3|\!|\!|\widetilde{p}_h|\!|\!|^2\leq0.
\end{equation*}
Therefore, we obtain $\widetilde{p}_h=0.$

Next, we denote by $\widetilde{T}_h$ the difference between two solutions of \eqref{s2}, then
\begin{equation*}
a_0(\frac{\widetilde{T}_h}{\tau},s_h)-(\mathcal M({\mathbf K}\nabla p_h^{n})\cdot\nabla \widetilde{T}_h,s_h)+c(\bm\Theta;\widetilde{T}_h,s_h)=0.
\end{equation*}
Choosing $s_h=\widetilde{T}_h$ in the equation above, together with Lemma \ref{lem0}, we get
\begin{align*}
&a_0(\frac{\widetilde{T}_h}{\tau},\widetilde{T}_h)+c(\bm\Theta;\widetilde{T}_h,\widetilde{T}_h)-(\mathcal M({\mathbf K}\nabla p_h^{n})\cdot\nabla \widetilde T_h, \widetilde T_h)\\
\geq &\frac{a_0}{\tau}||\widetilde{T}_h||^2+C_3|\!|\!|\widetilde{T}_h|\!|\!|^2-M||\nabla \widetilde T_h||||\widetilde T_h||\\
\geq& (\frac{a_0}{\tau}-\frac{M^2}{2C_3})||\widetilde{T}_h||^2+\frac{C_3}{2}|\!|\!|\widetilde{T}_h|\!|\!|^2.
\end{align*}
When $\frac{a_0}{\tau}-\frac{M^2}{2C_3}\ge0$, we have
\begin{equation*}
(\frac{a_0}{\tau}-\frac{M^2}{2C_3})||\widetilde{T}_h||^2+\frac{C_3}{2}|\!|\!|\widetilde{T}_h|\!|\!|^2\leq0.
\end{equation*}
Therefore, we obtain $\widetilde{T}_h=0.$

Finally, let $\widetilde{\mathbf u}_h$ be the difference between two solutions of \eqref{s3}, then
 \begin{equation*}
a(\widetilde{\mathbf u}_h,{\mathbf v_h} )+\gamma(\frac{\widetilde{\mathbf u}_h}{\tau},{\mathbf v_h})=0.
\end{equation*}
Choosing ${\mathbf v_h}=\widetilde{\mathbf u}_h$ in the equation above and applying Lemma \ref{lem2.5}, we have
\begin{equation*}
C_1|\!|\!|\widetilde{\mathbf u}_h|\!|\!|_V^2+\frac{\gamma}{\tau}||\widetilde{\mathbf u}_h||^2\leq0.
\end{equation*}
Therefore, we obtain $\widetilde{\mathbf u}_h=0.$ The uniqueness of $p_h^1$, $T_h^1$ and ${\mathbf u}_h^1$ can be similarly proved.
\end{proof}

Now, we are ready for the stability analysis of \eqref{s1}-\eqref{s3}. To this end, we introduce the following formula,
\begin{equation}
\sum\limits_{n=1}^{N-1}A^{n+1}(B^{n+1}-B^{n})=A^NB^N-A^1B^1-\sum\limits_{n=1}^{N-1}(A^{n+1}-A^{n})B^{n}.\label{f}
\end{equation}

\begin{theorem}\label{sta}
Assuming that (A1)-(A6) hold true. Let $\gamma$ be large enough, and $\gamma-\frac{\alpha^2+\beta^2}{C_3}(1+C)\ge0$. When $\tau\le\frac{C_3(a_0-b_0)}{4M^2}$, the system \eqref{s1}-\eqref{s3} is stable.
\end{theorem}
\begin{proof}
Taking ${\mathbf v}_h=\partial_\tau{\mathbf u}_h^{n+1}$, $q_h=p_h^{n+1}$ and $s_h=T_h^{n+1}$ in \eqref{s1},\eqref{s2} and \eqref{s3}, respectively, and adding these three equations yields that
\begin{align}
&\frac{c_0}{\tau}(p_h^{n+1}-p_h^{n},p_h^{n+1})+\frac{a_0}{\tau}(T_h^{n+1}-T_h^{n}, T_h^{n+1})+\frac{1}{\tau}a({\mathbf u}_h^{n+1},{\mathbf u}_h^{n+1}-{\mathbf u}_h^{n})\notag\\
& +c({\mathbf K};p_h^{n+1},p_h^{n+1})+c(\bm\Theta;T_h^{n+1},T_h^{n+1})+\gamma(\partial_\tau\mathbf u_h^{n+1},\partial_\tau\mathbf u_h^{n+1})-\gamma({\partial_\tau\mathbf u}_h^n,{\partial_\tau\mathbf u}_h^{n+1})\notag\\
=&\frac{b_0}{\tau}(T_h^{n}-T_h^{n-1}, p_h^{n+1})+\frac{b_0}{\tau}(p_h^{n+1}-p_h^{n},T_h^{n+1})+\alpha b(\partial_\tau{\mathbf u}_h^{n+1}-\partial_\tau{\mathbf u}_h^{n},p_h^{n+1})\notag\\
&+\beta b(\partial_\tau{\mathbf u}_h^{n+1}-\partial_\tau{\mathbf u}_h^{n},T_h^{n+1})+(\mathcal M({\mathbf K}\nabla p_h^n)\cdot\nabla  T_h^{n+1}, T_h^{n+1})\notag\\
&+(g^{n+1},p_h^{n+1})+(z^{n+1},T_h^{n+1})+(\bm f^{n+1},\partial_\tau\mathbf u_h^{n+1}).\label{sta1}
\end{align}
Note that
\begin{align*}
&\frac{c_0}{\tau}(p_h^{n+1}-p_h^{n},p_h^{n+1})=\frac{c_0}{2\tau}(||p_h^{n+1}||^2-||p_h^{n}||^2)+\frac{c_0\tau}{2}||\partial_\tau p_h^{n+1}||^2,\\
&\frac{a_0}{\tau}(T_h^{n+1}-T_h^{n},T_h^{n+1})=\frac{a_0}{2\tau}(||T_h^{n+1}||^2-||T_h^{n}||^2)+\frac{a_0\tau}{2}||\partial_\tau T_h^{n+1}||^2,\\
&a({\mathbf u}_h^{n+1},{\partial_\tau\mathbf u}_h^{n+1})=\frac{1}{2\tau}( a({\mathbf u}_h^{n+1},{\mathbf u}_h^{n+1})- a({\mathbf u}_h^{n},{\mathbf u}_h^{n}))+\frac{1}{2\tau} a({\mathbf u}_h^{n+1}-{\mathbf u}_h^{n},{\mathbf u}_h^{n+1}-{\mathbf u}_h^{n})\\
&\ \ \  \  \ \ \ \ \ \ \ \ \ \ \ \ \ \ \ \ \ \ge\frac{1}{2\tau}( a({\mathbf u}_h^{n+1},{\mathbf u}_h^{n+1})- a({\mathbf u}_h^{n},{\mathbf u}_h^{n}))+\frac{C_1\tau}{2}|\!|\!|\partial_\tau{\mathbf u}_h^{n+1}|\!|\!|_V^2,\\
 &c({\mathbf K};p_h^{n+1},p_h^{n+1})\ge C_3|\!|\!|p_h^{n+1} |\!|\!|^2,\\
 &c(\bm\Theta;T_h^{n+1},T_h^{n+1})\ge C_3|\!|\!|T_h^{n+1} |\!|\!|^2,\\
&\gamma({\partial_\tau\mathbf u}_h^{n+1},{\partial_\tau\mathbf u}_h^{n+1})-\gamma({\partial_\tau\mathbf u}_h^n,{\partial_\tau\mathbf u}_h^{n+1}) =\frac{\gamma}{2}(||{\partial_\tau\mathbf u}_h^{n+1}||^2-||{\partial_\tau\mathbf u}_h^{n}||^2)+\frac{\gamma}{2}||{\partial_\tau\mathbf u}_h^{n+1}-{\partial_\tau\mathbf u}_h^{n}||^2,
\end{align*}
and
\begin{align*}
&\frac{b_0}{\tau}(T_h^{n}-T_h^{n-1}, p_h^{n+1})+\frac{b_0}{\tau}(p_h^{n+1}-p_h^{n},T_h^{n+1})\\
=&b_0(\partial_\tau T_h^{n}-\partial_\tau T_h^{n+1},p_h^{n+1})+b_0(\partial_\tau T_h^{n+1},p_h^{n+1})+b_0(\partial_\tau P_h^{n+1},T_h^{n+1})\\
=&b_0(\partial_\tau T_h^{n}-\partial_\tau T_h^{n+1},p_h^{n+1})+\tau b_0(\partial_\tau T_h^{n+1},\partial_\tau p_h^{n+1})+b_0\partial_\tau( T_h^{n+1},p_h^{n+1}).
\end{align*}
Consequently,
\begin{align}
 &\frac{c_0}{2\tau}(||p_h^{n+1}||^2-||p_h^{n}||^2)+\frac{a_0}{2\tau}(||T_h^{n+1}||^2-||T_h^{n}||^2)+\frac{1}{2\tau}( a({\mathbf u}_h^{n+1},{\mathbf u}_h^{n+1})- a({\mathbf u}_h^{n},{\mathbf u}_h^{n}))\notag\\
&+\frac{c_0\tau}{2}||\partial_\tau p_h^{n+1}||^2 +\frac{a_0\tau}{2}||\partial_\tau T_h^{n+1}||^2+\frac{C_1\tau}{2}|\!|\!|\partial_\tau{\mathbf u}_h^{n+1}|\!|\!|_V^2+C_3|\!|\!|p_h^{n+1} |\!|\!|^2+C_3|\!|\!|T_h^{n+1} |\!|\!|^2\notag\\
&+\frac{\gamma}{2}(||{\partial_\tau\mathbf u}_h^{n+1}||^2-||{\partial_\tau\mathbf u}_h^{n}||^2)+\frac{\gamma}{2}||{\partial_\tau\mathbf u}_h^{n+1}-{\partial_\tau\mathbf u}_h^{n}||^2\notag\\
\le&b_0(\partial_\tau T_h^{n}-\partial_\tau T_h^{n+1},p_h^{n+1})+\tau b_0(\partial_\tau T_h^{n+1},\partial_\tau p_h^{n+1})+b_0\partial_\tau( T_h^{n+1},p_h^{n+1})\notag\\
&+\alpha b(\partial_\tau{\mathbf u}_h^{n+1}-\partial_\tau{\mathbf u}_h^{n},p_h^{n+1})+\beta b(\partial_\tau{\mathbf u}_h^{n+1}-\partial_\tau{\mathbf u}_h^{n},T_h^{n+1})\notag\\
&+(\mathcal M({\mathbf K}\nabla p_h^n)\cdot\nabla  T_h^{n+1}, T_h^{n+1})+(g^{n+1},p_h^{n+1})+(z^{n+1},T_h^{n+1})+(\bm f^{n+1},\partial_\tau\mathbf u_h^{n+1}).\label{sta2}
\end{align}
Applying the summation operator $\tau \sum\limits_{n=1}^{N-1}$ to the both sides of \eqref{sta2}, we have
\begin{align}
 &\frac{c_0}{2}(||p_h^{N}||^2-||p_h^{1}||^2)+\frac{a_0}{2}(||T_h^{N}||^2-||T_h^{1}||^2)+\frac{1}{2}( a({\mathbf u}_h^{N},{\mathbf u}_h^{N})- a({\mathbf u}_h^{1},{\mathbf u}_h^{1}))\notag\\
&+\tau \sum\limits_{n=1}^{N-1}[\frac{c_0\tau}{2}||\partial_\tau p_h^{n+1}||^2 +\frac{a_0\tau}{2}||\partial_\tau T_h^{n+1}||^2+\frac{C_1\tau}{2}|\!|\!|\partial_\tau{\mathbf u}_h^{n+1}|\!|\!|_V^2+C_3|\!|\!|p_h^{n+1} |\!|\!|^2\notag\\
&+C_3|\!|\!|T_h^{n+1} |\!|\!|^2+\frac{\gamma}{2}||{\partial_\tau\mathbf u}_h^{n+1}-{\partial_\tau\mathbf u}_h^{n}||^2]+\frac{\gamma}{2}\tau(||{\partial_\tau\mathbf u}_h^{N}||^2-||{\partial_\tau\mathbf u}_h^{1}||^2)\notag\\
\le&b_0( T_h^{N},p_h^{N})-b_0( T_h^{1},p_h^{1})+ b_0\tau^2 \sum\limits_{n=1}^{N-1}(\partial_\tau T_h^{n+1},\partial_\tau p_h^{n+1})\notag\\
&+\alpha\tau \sum\limits_{n=1}^{N-1} b(\partial_\tau{\mathbf u}_h^{n+1}-\partial_\tau{\mathbf u}_h^{n},p_h^{n+1})+\beta\tau \sum\limits_{n=1}^{N-1} b(\partial_\tau{\mathbf u}_h^{n+1}-\partial_\tau{\mathbf u}_h^{n},T_h^{n+1})\notag\\
&+\tau \sum\limits_{n=1}^{N-1}(\mathcal M({\mathbf K}\nabla p_h^{n})\cdot\nabla  T_h^{n+1}, T_h^{n+1})+\tau \sum\limits_{n=1}^{N-1}[(g^{n+1},p_h^{n+1})+(z^{n+1},T_h^{n+1})]\notag\\
&+\tau \sum\limits_{n=1}^{N-1}(\bm f^{n+1},\partial_\tau\mathbf u_h^{n+1})+b_0\tau \sum\limits_{n=1}^{N-1}(\partial_\tau T_h^{n}-\partial_\tau T_h^{n+1},p_h^{n+1}).\label{sta3}
  \end{align}
Applying Cauchy-Schwartz inequality and Young's inequality, let us estimate the right side of \eqref{sta3} term by term,
\begin{align}
b_0( T_h^{N},p_h^{N})\le& \frac{b_0}{2}||T_h^{N}||^2+\frac{b_0}{2}||p_h^{N}||^2,\label{sdw1}\\
b_0\tau^2 \sum\limits_{n=1}^{N-1}(\partial_\tau T_h^{n+1},\partial_\tau p_h^{n+1})\le&\tau \sum\limits_{n=1}^{N-1}\frac{b_0}{2}\tau||\partial_\tau T_h^{n+1}||^2+\tau \sum\limits_{n=1}^{N-1}\frac{b_0}{2}\tau||\partial_\tau p_h^{n+1}||^2,\label{sdw2}\\
\tau \sum\limits_{n=1}^{N-1}(\mathcal M({\mathbf K}\nabla p_h^{n})\cdot\nabla  T_h^{n+1}, T_h^{n+1})
\le& \tau \sum\limits_{n=1}^{N-1}(\frac{C_3}{4}|\!|\!|T_h^{n+1}|\!|\!|^2+\frac{M^2}{C_3}||T_h^{n+1}||^2),\label{sdw5}\\
\tau \sum\limits_{n=1}^{N-1}[(g^{n+1},p_h^{n+1})+(z^{n+1},T_h^{n+1})]
\le& C\tau \sum\limits_{n=1}^{N-1}\bigg(||g^{n+1}||^2+||z^{n+1}||^2\bigg)+\frac{c_0-2b_0}{2}\tau||p_h^{N}||^2\notag\\
&+\frac{M^2}{C_3}\tau||T_h^{N}||^2+C\tau \sum\limits_{n=1}^{N-2}\bigg(||p_h^{n+1}||^2+||T_h^{n+1}||^2\bigg).\label{sdw6}
\end{align}
Since
\begin{align*}
&\alpha b(\partial_\tau{\mathbf u}_h^{n+1}-\partial_\tau{\mathbf u}_h^{n},p_h^{n+1})\\
=&\alpha\sum\limits_{K\in \mathcal T}\int_K (\partial_\tau{\mathbf u}_h^{n+1}-\partial_\tau{\mathbf u}_h^{n})\cdot\nabla p_h^{n+1}dK+\alpha\sum\limits_{e\in \Gamma}\int_{e}\{(\partial_\tau{\mathbf u}_h^{n+1}-\partial_\tau{\mathbf u}_h^{n})\}\cdot{\bm n_e}[p_h^{n+1}]ds\\
\le& \frac{C_3}{2}||\nabla p_h^{n+1} ||^2+\frac{\alpha^2}{2C_3}||\partial_\tau{\mathbf u}_h^{n+1}-\partial_\tau{\mathbf u}_h^{n}||^2+\frac{\alpha^2}{2C_3}\sum\limits_{e\in \Gamma}h_e||\partial_\tau{\mathbf u}_h^{n+1}-\partial_\tau{\mathbf u}_h^{n}||_{0,e}^2\\
&+\frac{C_3}{2}\sum\limits_{e\in \Gamma}h_e^{-1}||[p_h^{n+1}]||_{0,e}^2\\
\le& \frac{C_3}{2}|\!|\!|p_h^{n+1} |\!|\!|^2+\frac{\alpha^2}{2C_3}||\partial_\tau{\mathbf u}_h^{n+1}-\partial_\tau{\mathbf u}_h^{n}||^2+\frac{\alpha^2}{2C_3}Ch^2||\partial_\tau{\nabla\mathbf u}_h^{n+1}-\partial_\tau{\nabla\mathbf u}_h^{n}||^2 \\
\le& \frac{C_3}{2}|\!|\!|p_h^{n+1} |\!|\!|^2+\frac{\alpha^2}{2C_3}(1+C)||\partial_\tau{\mathbf u}_h^{n+1}-\partial_\tau{\mathbf u}_h^{n}||^2,
\end{align*}
we provide
\begin{align}
&\alpha\tau \sum\limits_{n=1}^{N-1} b(\partial_\tau{\mathbf u}_h^{n+1}-\partial_\tau{\mathbf u}_h^{n},p_h^{n+1})+\beta\tau \sum\limits_{n=1}^{N-1} b(\partial_\tau{\mathbf u}_h^{n+1}-\partial_\tau{\mathbf u}_h^{n},T_h^{n+1})\\
\le&\frac{\alpha^2+\beta^2}{2C_3}(1+C)\tau \sum\limits_{n=1}^{N-1}||\partial_\tau{\mathbf u}_h^{n+1}-\partial_\tau{\mathbf u}_h^{n}||^2+\frac{C_3}{2}\tau \sum\limits_{n=1}^{N-1}|\!|\!|p_h^{n+1}|\!|\!|^2+\frac{C_3}{2}\tau \sum\limits_{n=1}^{N-1}|\!|\!|T_h^{n+1}|\!|\!|^2.
\end{align}
Applying \eqref{f}, we obtain
\begin{align}
&\tau \sum\limits_{n=1}^{N-1}(\bm f^{n+1},\partial_\tau\mathbf u_h^{n+1})\notag\\
=&(\bm f^{N},\mathbf u_h^{N})-(\bm f^{1},\mathbf u_h^{1})-\tau \sum\limits_{n=1}^{N-1}(\partial_\tau\bm f^{n+1},\mathbf u_h^{n})\notag\\
\le& \frac{1}{C_1}||\bm f^{N}||^2+\frac{C_1}{4}|\!|\!|\mathbf u_h^{N}|\!|\!|_V^2-(\bm f^{1},\mathbf u_h^{1})+C\tau \sum\limits_{n=1}^{N-1}||\partial_\tau\bm f^{n+1}||^2+C\tau \sum\limits_{n=1}^{N-1}||\mathbf u_h^{n}||^2.
\end{align}
Similarly,
\begin{align}
&b_0\tau \sum\limits_{n=1}^{N-1}(\partial_\tau T_h^{n}-\partial_\tau T_h^{n+1},p_h^{n+1})\notag\\
=&b_0\tau(\partial_\tau T_h^{1},p_h^{1})-b_0\tau(\partial_\tau T_h^{N},p_h^{N})+b_0\tau  \sum\limits_{n=1}^{N-1}\tau(\partial_\tau T_h^{n},\partial_\tau p_h^{n+1})\notag\\
\le&\frac{b_0}{2}\tau^2||\partial_\tau T_h^{1}||^2+\frac{b_0}{2}||p_h^{1}||^2+\frac{b_0}{2}\tau^2||\partial_\tau T_h^{N}||^2+\frac{b_0}{2}||p_h^{N}||^2\notag\\
&+\frac{b_0}{2}\tau \sum\limits_{n=1}^{N-1}\tau||\partial_\tau T_h^{n}||^2+\frac{b_0}{2}\tau \sum\limits_{n=1}^{N-1}\tau||\partial_\tau p_h^{n+1}||^2.\label{sta5}
\end{align}
Substituting \eqref{sdw1}-\eqref{sta5} into \eqref{sta3}, we have
\begin{align*}
 &\frac{c_0-2b_0}{2}||p_h^{N}||^2+\frac{a_0-b_0}{2}||T_h^{N}||^2+\frac{C_1}{4}|\!|\!|{\mathbf u}_h^{N}|\!|\!|_V^2+\frac{\gamma}{2}\tau||{\partial_\tau\mathbf u}_h^{N}||^2\notag\\
&+\tau \sum\limits_{n=1}^{N-1}\bigg[\frac{(c_0-2b_0)\tau}{2}||\partial_\tau p_h^{n+1}||^2 +\frac{(a_0-2b_0)\tau}{2}||\partial_\tau T_h^{n+1}||^2+\frac{C_1\tau}{2}|\!|\!|\partial_\tau{\mathbf u}_h^{n+1}|\!|\!|_V^2\notag\\
&+\frac{C_3}{2}|\!|\!|p_h^{n+1} |\!|\!|^2+\frac{C_3}{4}|\!|\!|T_h^{n+1} |\!|\!|^2+\frac{1}{2}\left(\gamma-\frac{\alpha^2+\beta^2}{C_3}(1+C)\right)||{\partial_\tau\mathbf u}_h^{n+1}-{\partial_\tau\mathbf u}_h^{n}||^2\bigg]\notag\\
\le&\frac{c_0+b_0}{2}||p_h^{1}||^2+\frac{a_0}{2}||T_h^{1}||^2+\frac{b_0}{2}\tau^2||\partial_\tau T_h^{1}||^2+C|\!|\!|{\mathbf u}_h^{1}|\!|\!|_V^2-b_0( T_h^{1},p_h^{1})-(\bm f^{1},\mathbf u_h^{1})\notag\\
&+\frac{\gamma}{2}\tau||{\partial_\tau\mathbf u}_h^{1}||^2+C||\bm f^{N}||^2+C\tau \sum\limits_{n=1}^{N-1}(||\partial_\tau\bm f^{n+1}||^2+||g^{n+1}||^2+||z^{n+1}||^2)\notag\\
&+C\tau \sum\limits_{n=1}^{N-1}(||\mathbf u_h^{n}||^2+||p_h^{n}||^2+||T_h^{n}||^2)+\frac{c_0-2b_0}{2}\tau||p_h^{N}||^2+\frac{2M^2}{C_3}\tau||T_h^{N}||^2.
\end{align*}
Choosing $\gamma$ large enough, when $c_0-2b_0>0$, $a_0-2b_0\ge0$, $\gamma-\frac{\alpha^2+\beta^2}{C_3}(1+C)\ge0$ and $\tau\le\frac{C_3(a_0-b_0)}{4M^2}$, by Gronwall inequality, we provide
\begin{align*}
&\frac{c_0-2b_0}{2}(1-\tau)||p_h^{N}||^2+\left(\frac{a_0-b_0}{2}-\frac{2M^2}{C_3}\tau\right)||T_h^{N}||^2+\frac{C_1}{4}|\!|\!|{\mathbf u}_h^{N}|\!|\!|_V^2+\frac{\gamma}{2}\tau||{\partial_\tau\mathbf u}_h^{N}||^2\notag\\
\le&\frac{c_0}{2}||p_h^{1}||^2+\frac{a_0}{2}||T_h^{1}||^2+\frac{b_0}{2}||T_h^{0}||^2+C|\!|\!|{\mathbf u}_h^{1}|\!|\!|_V^2+\frac{\gamma}{2}\tau||{\partial_\tau\mathbf u}_h^{1}||^2\notag\\
&+C||\bm f^{N}||^2-C||\bm f^{1}||^2+C\tau \sum\limits_{n=1}^{N-1}(||\partial_\tau\bm f^{n+1}||^2+||g^{n+1}||^2+||z^{n+1}||^2).
\end{align*}
The proof is completed.
\end{proof}

\section{Error estimate}\label{EA}
The following error estimation will make use of several projection operators and lemmas that we describe in this section. Then, we build the related error equations and obtain optimal convergence order estimates for our proposed SDG approach .

\subsection{Some projection operators}\label{EA1}
Let  \begin{align*}
&{\mathbf u}(t)-{\mathbf u}_h(t)=[{\mathbf u}(t)-R_u{\mathbf u}(t)]+[R_u{\mathbf u}(t)-{\mathbf u}_h(t)]={\bm\rho}_u(t)+{\bm\theta}_u(t),\\
&p(t)-p_h(t)=[p(t)-R_pp(t)]+[R_pp(t)-p_h(t)]=\rho_p(t)+\theta_p(t),\\
&T(t)-T_h(t)=[T(t)-R_TT(t)]+[R_TT(t)-T_h(t)]=\rho_T(t)+\theta_T(t),
\end{align*}
where the projection operators $R_u$, $R_p$ and $R_T$ are defined in \eqref{ea1}-\eqref{ea3}.

Next, the following lemmas about the approximation capability of elliptic projections are stated.
\begin{lemma}\label{lem3.7}\cite{25} If we choose $\sigma_1$ and $\sigma_2$ large enough, the following estimates hold for ${\mathbf u}\in (H^{k_1+1}(\Omega))^d$, $p\in H^{k_2+1}(\Omega)$ and $T\in H^{k_3+1}(\Omega)$,
\begin{align*}
\|{\bm\rho}_u\|+h|\!|\!|{\bm\rho}_u|\!|\!|_V&\leq Ch^{k_1+1}\|{\mathbf u}\|_{k_1+1},\\
\|{\rho}_p\|+h|\!|\!|{\rho}_p|\!|\!|&\leq Ch^{k_2+1}\|p\|_{k_2+1},\\
\|{\rho}_T\|+h|\!|\!|{\rho}_T|\!|\!|&\leq Ch^{k_3+1}\|T\|_{k_3+1}.
\end{align*}
\end{lemma}

\begin{lemma}\label{lem3.8}\cite{MR3742883} For ${\mathbf u}\in (H^{2}(\Omega))^d$, $p\in H^{2}(\Omega)$ and $T\in H^{2}(\Omega)$, the following inequalities hold
\begin{align*}
\|R_u{\mathbf u}\|&\leq C\|{\mathbf u}\|_{2},\\
\|R_pp\|&\leq C\|p\|_{2},\\
\|R_TT\|&\leq C\|T\|_{2}.
\end{align*}
\end{lemma}

Also, it should be noted that the following results hold,
\begin{align}
\|\partial_tp^{n+1}-\partial_\tau p^{n+1}\|^2
&\leq C\tau\int_{t_{n}}^{t_{n+1}}\|p_{tt}\|^2dt,\label{pt}\\
||p^{n+1}-p^{n}||^2&\le C\tau\int_{t_n}^{t_{n+1}}||p_t||^2dt,\label{pt0}\\
||\partial_t{\mathbf u}^{n+1}-\partial_\tau{\mathbf u}^{n+1}||^2&\le C\tau\int_{t_n}^{t_{n+1}}||{\mathbf u}_{tt}||^2dt,\\
||{\mathbf u}^{n+1}-{\mathbf u}^{n}||^2&\le C\tau\int_{t_n}^{t_{n+1}}||{\mathbf u}_{t}||^2dt,
\end{align}
where $\partial_\tau p^{n+1}=(p^{n+1}-p^{n})/\tau$ and $\partial_\tau{\mathbf u}^{n+1}=({\mathbf u}^{n+1}-{\mathbf u}^{n})/\tau$.

\subsection{Optimal convergence order estimates}
In this subsection, we will derive the optimal error estimates for the SDG method \eqref{s1}-\eqref{s3}.
\begin{lemma}\label{lemxr}\cite{chen} (Consistency) If ${\mathbf u}$, $p$ and $T$ is the solution of \eqref{e1}-\eqref{e7} and essentially bounded, for any $q_h$, $s_h\in H^{2}(\mathcal T_h)$ and ${\mathbf v}_h\in (H^{2}(\mathcal T_h))^d,$ there hold
\begin{align}
&c_0(p_t^{n+1},q_h)-b_0(T_t^{n+1},q_h)+\alpha b({\mathbf u}_t^{n+1},q_h)+c({\mathbf K};p^{n+1},q_h)=(g^{n+1},q_h),\label{b1}\\
 & a_0(T_t^{n+1},s_h)-b_0(p_t^{n+1},s_h)+\beta b({\mathbf u}_t^{n+1},s_h)+c(\bm\Theta;T^{n+1},s_h)\notag\\
  &\hspace{5cm}-({\mathbf K}\nabla p^{n+1}\cdot\nabla T^{n+1},s_h)=(z^{n+1},s_h), \label{b2}\\
 & a({\mathbf u}^{n+1},{\mathbf v}_h )-\alpha b({\mathbf v}_h,p^{n+1})-\beta b({\mathbf v}_h,T^{n+1}) = ({\bm f}^{n+1},{\mathbf v}_h). \label{b3}
\end{align}
\end{lemma}

Subtracting \eqref{s1} from \eqref{b1}, \eqref{s2} from \eqref{b2} and \eqref{s3} from \eqref{b3}, together with \eqref{ea1}, \eqref{ea2} and \eqref{ea3}, we obtain the error equations
\begin{align}
&c_0(\partial_\tau \theta_p^{n+1},q_h)-b_0(\partial_\tau \theta_T^n,q_h)+\alpha b(\partial_\tau{\bm \theta}_u^n,q_h)+c({\mathbf K};\theta_p^{n+1},q_h)\notag\\
=&-c_0(p_t^{n+1}-\partial_\tau p^{n+1},q_h)-c_0(\partial_\tau\rho_p^{n+1},q_h)\notag\\
&+b_0(T_t^{n+1}-T_t^n,q_h)+b_0(T_t^n-\partial_\tau T^n,q_h)+b_0(\partial_\tau\rho_T^n,q_h)\notag\\
&-\alpha b({\mathbf u}_t^{n+1}-{\mathbf u}_t^n,q_h)-\alpha b({\mathbf u}_t^n-\partial_\tau {\mathbf u}^n,q_h)-\alpha b(\partial_\tau\bm\rho_u^n,q_h),\label{ea4}\\
&a_0(\partial_\tau \theta_T^{n+1},s_h)-b_0(\partial_\tau \theta_p^{n+1},s_h)+\beta b(\partial_\tau{\bm \theta}_u^n,s_h)+c(\bm\Theta;\theta_T^{n+1},s_h)\notag\\
=&-a_0(T_t^{n+1}-\partial_\tau T^{n+1},s_h)-a_0(\partial_\tau\rho_T^{n+1},s_h)+b_0(p_t^{n+1}-\partial_\tau p^{n+1},s_h)+b_0(\partial_\tau\rho_p^{n+1},s_h)\notag\\
&-\beta b({\mathbf u}_t^{n+1}-{\mathbf u}_t^n,s_h)-\beta b({\mathbf u}_t^{n}-\partial_\tau {\mathbf u}^n,s_h)-\beta b(\partial_\tau\bm\rho_u^n,s_h)\notag\\
&+(\mathcal M({\mathbf K}\nabla p_h^{n})\cdot\nabla \theta_T^{n+1},s_h)+(({\mathbf K}\nabla p^{n+1}-\mathcal M({\mathbf K}\nabla p_h^{n}))\cdot\nabla R_TT^{n+1},s_h),\label{ea5}\\
&a({\bm \theta}_u^{n+1},{\mathbf v_h} )-\alpha b({\mathbf v_h},\theta_p^{n+1})-\beta b({\mathbf v_h},\theta_T^{n+1})+\gamma(\partial_\tau{\bm\theta}_u^{n+1},{\mathbf v_h})-\gamma(\partial_\tau{\bm\theta}_u^{n},{\mathbf v_h})\notag\\
=&\alpha b({\mathbf v_h},\rho_p^{n+1})+\beta b({\mathbf v_h},\rho_T^{n+1})+\gamma(\partial_\tau R_u{\mathbf u}^{n+1}-\partial_\tau R_u{\mathbf u}^{n},{\mathbf v_h}).\label{ea6}
\end{align}

Now we are ready for the optimal convergence order estimates as follows.
\begin{theorem}\label{thm1}
For $k_1$, $k_2$, $k_3\ge1,$ let $({\mathbf u},p,T) \in L^{\infty}(0,t_f;(H^{k_1+1}(\Omega)\cap H_0^1(\Omega))^d)\times L^{\infty}(0,t_f;\\
H^{k_2+1}(\Omega)\cap H_0^1(\Omega))\times L^{\infty}(0,t_f;H^{k_3+1}(\Omega)\cap H_0^1(\Omega))$ and $({\mathbf u}_h^n,p_h^n,T_h^n)\in {\mathbf V}_h\times Q_h\times S_h$ be the solutions of (\ref{b1})-(\ref{b3}) and (\ref{s1})-(\ref{s3}), respectively. Assuming that $\gamma$, $\sigma_1$, $\sigma_2$ are large enough, then for $n\geq2$,
\begin{align*}
&|\!|\!|\mathbf u^n-\mathbf u_h^n|\!|\!|_V+||p^{n}-p_h^n||+||T^{n}-T_h^n||\\
\le& C(\tau+h^{k_1}+h^{k_2+1}+h^{k_3+1})\\
&+(C||p_h^{1}-R_pp^1||^2+C||T_h^{1}-R_TT^1||^2+C|\!|\!|\bm u_h^{1}-R_u\mathbf u^1|\!|\!|_V^2+\gamma\tau||\frac{\mathbf u_h^{1}-R_u\mathbf u^1}{\tau_0}||^2)^{\frac{1}{2}}.
\end{align*}
\end{theorem}
\begin{proof}
\indent Taking ${\mathbf v}_h=\partial_\tau{\bm\theta}_u^{n+1}$, $q_h=\theta_p^{n+1}$ and $s_h=\theta_T^{n+1}$ in \eqref{ea4}, \eqref{ea5} and \eqref{ea6}, respectively, and adding these three equations yields that
\begin{align}
 &c_0(\partial_\tau \theta_p^{n+1},\theta_p^{n+1})+a_0(\partial_\tau \theta_T^{n+1}, \theta_T^{n+1})+a({\bm\theta}_u^{n+1},{\partial_\tau\bm\theta}_u^{n+1})\notag\\
&+c({\mathbf K};\theta_p^{n+1},\theta_p^{n+1}) +c(\bm\Theta;\theta_T^{n+1},\theta_T^{n+1})+\gamma({\partial_\tau\bm\theta}_u^{n+1}-{\partial_\tau\bm\theta}_u^n,{\partial_\tau\bm\theta}_u^{n+1})\notag\\
 =  &b_0(\partial_\tau \theta_T^{n}, \theta_p^{n+1})+b_0(\partial_\tau\theta_p^{n+1}, \theta_T^{n+1})-c_0(p_t^{n+1}-\partial_\tau p^{n+1}, \theta_p^{n+1})-c_0(\partial_\tau\rho_p^{n+1},\theta_p^{n+1})\notag\\
 &+b_0(T_t^{n+1}-T_t^n, \theta_p^{n+1})+b_0(T_t^{n}-\partial_\tau T^{n}, \theta_p^{n+1})+b_0(\partial_\tau\rho_T^n,\theta_p^{n+1})\notag\\
 &-a_0(T_t^{n+1}-\partial_\tau T^{n+1},\theta_T^{n+1})-a_0(\partial_\tau\rho_T^{n+1},\theta_T^{n+1})\notag\\
 &+b_0(p_t^{n+1}-\partial_\tau p^{n+1},\theta_T^{n+1})+b_0(\partial_\tau \rho_p^{n+1},\theta_T^{n+1})\notag\\
 &-\alpha b({\mathbf u}_t^{n+1}-{\mathbf u}_t^n,\theta_p^{n+1})-\alpha b({\mathbf u}_t^{n}-\partial_\tau {\mathbf u}^{n},\theta_p^{n+1})-\alpha b(\partial_\tau\bm\rho_u^{n},\theta_p^{n+1})\notag\\
 &-\beta b({\mathbf u}_t^{n+1}-{\mathbf u}_t^n,\theta_T^{n+1})-\beta b({\mathbf u}_t^n-\partial_\tau {\mathbf u}^n,\theta_T^{n+1})-\beta b(\partial_\tau\bm\rho_u^n,\theta_T^{n+1})\notag\\
 &+(\mathcal M({\mathbf K}\nabla p_h^{n})\cdot\nabla \theta_T^{n+1},\theta_T^{n+1})+(({\mathbf K}\nabla p^{n+1}-\mathcal M({\mathbf K}\nabla p_h^{n}))\cdot\nabla R_TT^{n+1},\theta_T^{n+1})\notag\\
  &+\alpha(b(\partial_\tau\bm\theta_u^{n+1},\theta_p^{n+1})-b(\partial_\tau\bm\theta_u^n,\theta_p^{n+1}))+\beta(b(\partial_\tau\bm\theta_u^{n+1},\theta_T^{n+1})-b(\partial_\tau\bm\theta_u^n,\theta_T^{n+1}))\notag\\
 &+\alpha b(\partial_\tau\bm\theta_u^{n+1},\rho_p^{n+1})+\beta b(\partial_\tau\bm\theta_u^{n+1},\rho_T^{n+1})+\gamma(\partial_\tau R_u{\mathbf u}^{n+1}-\partial_\tau R_u{\mathbf u}^{n},\partial_\tau\bm\theta_u^{n+1}).\label{ea7}
  \end{align}
Applying the symmetry of $a(\cdot,\cdot)$ and the coercivity results in Lemma \eqref{lem2.4}, we get
\begin{align*}
 &c_0(\partial_\tau \theta_p^{n+1},\theta_p^{n+1})=\frac{c_0}{2\tau}(||\theta_p^{n+1}||^2-||\theta_p^{n}||^2)+\frac{c_0}{2\tau}||\theta_p^{n+1}-\theta_p^{n}||^2,\\
& a_0(\partial_\tau \theta_T^{n+1}, \theta_T^{n+1})=\frac{a_0}{2\tau}(||\theta_T^{n+1}||^2-||\theta_T^{n}||^2)+\frac{a_0}{2\tau}||\theta_T^{n+1}-\theta_T^{n}||^2,\\
& a({\bm\theta}_u^{n+1},{\partial_\tau\bm\theta}_u^{n+1})=\frac{1}{2\tau}( a({\bm\theta}_u^{n+1},{\bm\theta}_u^{n+1})-a({\bm\theta}_u^{n},{\bm\theta}_u^{n}))+\frac{1}{2\tau} a({\bm\theta}_u^{n+1}-{\bm\theta}_u^{n},{\bm\theta}_u^{n+1}-{\bm\theta}_u^{n})\\
&\hspace{2.5cm} \ge\frac{1}{2\tau}( a({\bm\theta}_u^{n+1},{\bm\theta}_u^{n+1})-  a({\bm\theta}_u^{n},{\bm\theta}_u^{n}))+\frac{C_1\tau}{2}|\!|\!|\partial_\tau{\bm\theta}_u^{n+1}|\!|\!|_V^2,\\
 &c({\mathbf K};\theta_p^{n+1},\theta_p^{n+1})\ge C_3|\!|\!|\theta_p^{n+1} |\!|\!|^2,\\
& c(\bm\Theta;\theta_T^{n+1},\theta_T^{n+1})\ge C_3|\!|\!|\theta_T^{n+1} |\!|\!|^2,\\
&\gamma({\partial_\tau\bm\theta}_u^{n+1},{\partial_\tau\bm\theta}_u^{n+1})-\gamma({\partial_\tau\bm\theta}_u^n,{\partial_\tau\bm\theta}_u^{n+1}) =\frac{\gamma}{2}(||{\partial_\tau\bm\theta}_u^{n+1}||^2-||{\partial_\tau\bm\theta}_u^{n}||^2)+\frac{\gamma}{2}||{\partial_\tau\bm\theta}_u^{n+1}-{\partial_\tau\bm\theta}_u^{n}||^2.
\end{align*}
The first item on the right side of \eqref{ea7} can be bounded as follows,
\begin{align*}
&b_0(\partial_\tau \theta_T^n,\theta_p^{n+1})+b_0(\partial_\tau\theta_p^{n+1},\theta_T^{n+1})\\
=&b_0(\partial_\tau \theta_T^{n+1},\theta_p^{n+1})+b_0(\partial_\tau\theta_p^{n+1},\theta_T^{n+1})+b_0(\partial_\tau \theta_T^{n}-\partial_\tau \theta_T^{n+1},\theta_p^{n+1})\\
=&b_0\tau(\partial_\tau \theta_T^{n+1},\partial_\tau\theta_p^{n+1})+b_0\partial_\tau(\theta_p^{n+1},\theta_T^{n+1})+b_0(\partial_\tau \theta_T^{n}-\partial_\tau \theta_T^{n+1},\theta_p^{n+1})\\
\le&\frac{b_0}{2}\tau||\partial_\tau \theta_T^{n+1}||^2+\frac{b_0}{2}\tau||\partial_\tau \theta_p^{n+1}||^2+b_0\partial_\tau(\theta_p^{n+1},\theta_T^{n+1})+b_0(\partial_\tau \theta_T^{n}-\partial_\tau \theta_T^{n+1},\theta_p^{n+1}).
\end{align*}
Substituting the above equations into \eqref{ea7} and applying the summation operator $2\tau \sum\limits_{n=1}^{N-1}$ to the both sides of the obtained inequality, we obtain
\begin{align}\label{z11}
&c_0(||\theta_p^{N}||^2-||\theta_p^{1}||^2)+(c_0-b_0)\tau^2\sum\limits_{n=1}^{N-1}||\partial_\tau\theta_p^{n+1}||^2+a_0(||\theta_T^{N}||^2-||\theta_T^{1}||^2)\notag\\
&+(a_0-b_0)\tau^2\sum\limits_{n=1}^{N-1}||\partial_\tau\theta_T^{n+1}||^2+a({\bm\theta}_u^{N},{\bm\theta}_u^{N})-a({\bm\theta}_u^{1},{\bm\theta}_u^{1})+C_1\tau^2\sum\limits_{n=1}^{N-1}|\!|\!|\partial_\tau{\bm\theta}_u^{n+1}|\!|\!|_V^2 \notag\\
&+2C_3\tau\sum\limits_{n=1}^{N-1}(|\!|\!|\theta_p^{n+1} |\!|\!|^2+|\!|\!|\theta_T^{n+1} |\!|\!|^2)+\gamma\tau(||{\partial_\tau\bm\theta}_u^{N}||^2-||{\partial_\tau\bm\theta}_u^{1}||^2)+\gamma\tau\sum\limits_{n=1}^{N-1}||{\partial_\tau\bm\theta}_u^{n+1}-{\partial_\tau\bm\theta}_u^{n}||^2\notag\\
\le&2b_0((\theta_p^{N}, \theta_T^{N})-(\theta_p^{1}, \theta_T^{1}))+2b_0\tau\sum\limits_{n=1}^{N-1}(\partial_\tau\theta_T^{n}-\partial_\tau\theta_T^{n+1}, \theta_p^{n+1})\notag\\
 &-2\tau\sum\limits_{n=1}^{N-1}[c_0(p_t^{n+1}-\partial_\tau p^{n+1}, \theta_p^{n+1})+c_0(\partial_\tau\rho_p^{n+1},\theta_p^{n+1})]\notag\\
 &+2\tau\sum\limits_{n=1}^{N-1}[b_0(T_t^{n+1}-T_t^n, \theta_p^{n+1})+b_0(T_t^{n}-\partial_\tau T^{n}, \theta_p^{n+1})+b_0(\partial_\tau\rho_T^n,\theta_p^{n+1})]\notag\\
 &-2\tau\sum\limits_{n=1}^{N-1}[a_0(T_t^{n+1}-\partial_\tau T^{n+1},\theta_T^{n+1})+a_0(\partial_\tau\rho_T^{n+1},\theta_T^{n+1})]\notag\\
 &+2\tau\sum\limits_{n=1}^{N-1}[b_0(p_t^{n+1}-\partial_\tau p^{n+1},\theta_T^{n+1})+b_0(\partial_\tau \rho_p^{n+1},\theta_T^{n+1})]\notag\\
 &-2\tau\sum\limits_{n=1}^{N-1}[\alpha b({\mathbf u}_t^{n+1}-{\mathbf u}_t^n,\theta_p^{n+1})+\alpha b({\mathbf u}_t^{n}-\partial_\tau {\mathbf u}^{n},\theta_p^{n+1})+\alpha b(\partial_\tau\bm\rho_u^{n},\theta_p^{n+1})]\notag\\
 &- 2\tau\sum\limits_{n=1}^{N-1}[\beta b({\mathbf u}_t^{n+1}-{\mathbf u}_t^n,\theta_T^{n+1})+\beta b({\mathbf u}_t^n-\partial_\tau {\mathbf u}^n,\theta_T^{n+1})+\beta b(\partial_\tau\bm\rho_u^n,\theta_T^{n+1})]\notag\\
 &+2\tau\sum\limits_{n=1}^{N-1}(\mathcal M({\mathbf K}\nabla p_h^{n})\cdot\nabla \theta_T^{n+1},\theta_T^{n+1})+2\tau\sum\limits_{n=1}^{N-1}(({\mathbf K}\nabla p^{n+1}-\mathcal M({\mathbf K}\nabla p_h^{n}))\cdot\nabla R_TT^{n+1},\theta_T^{n+1})\notag\\
  &+2\tau\sum\limits_{n=1}^{N-1}\alpha(b(\partial_\tau\bm\theta_u^{n+1},\theta_p^{n+1})-b(\partial_\tau\bm\theta_u^n,\theta_p^{n+1}))+2\tau\sum\limits_{n=1}^{N-1}\beta(b(\partial_\tau\bm\theta_u^{n+1},\theta_T^{n+1})-b(\partial_\tau\bm\theta_u^n,\theta_T^{n+1}))\notag\\
 &+2\tau\sum\limits_{n=1}^{N-1}\alpha b(\partial_\tau\bm\theta_u^{n+1},\rho_p^{n+1})+2\tau\sum\limits_{n=1}^{N-1}\beta b(\partial_\tau\bm\theta_u^{n+1},\rho_T^{n+1})\notag\\
 &+2\tau\sum\limits_{n=1}^{N-1}\gamma(\partial_\tau R_u{\mathbf u}^{n+1}-\partial_\tau R_u{\mathbf u}^{n},\partial_\tau\bm\theta_u^{n+1})\notag\\
  :=&I_1+\cdots+I_{15}.
\end{align}

Next, we proceed by finding upper bounds for the terms $I_1$-$I_{15}$ on the right side of \eqref{z11}. Applying Cauchy-Schwarz inequality and Young's inequality, we have
\begin{align*}
|I_1|=2b_0|(\theta_p^{N},\theta_T^{N})-(\theta_p^{1},\theta_T^{1})|\le b_0||\theta_p^{N}||^2+b_0||\theta_T^{N}||^2+b_0||\theta_p^{1}||^2+b_0||\theta_T^{1}||^2.
\end{align*}
Applying \eqref{f}, we arrive at
\begin{align*}
|I_2|=&2b_0\tau|\sum\limits_{n=1}^{N-1}(\partial_\tau \theta_T^{n}-\partial_\tau \theta_T^{n+1},\theta_p^{n+1})|\\
=&2b_0\tau|(\partial_\tau \theta_T^{1},\theta_p^{1})-(\partial_\tau \theta_T^{N},\theta_p^{N})+\sum\limits_{n=1}^{N-1}\tau(\partial_\tau \theta_T^{n},\partial_\tau \theta_p^{n+1})|\\
\le&b_0\tau^2||\partial_\tau \theta_T^{1}||^2+b_0||\theta_p^{1}||^2+b_0\tau^2||\partial_\tau \theta_T^{N}||^2+b_0||\theta_p^{N}||^2\\
&+b_0\tau \sum\limits_{n=1}^{N-1}\tau||\partial_\tau \theta_T^{n}||^2+b_0\tau \sum\limits_{n=1}^{N-1}\tau||\partial_\tau \theta_p^{n+1}||^2.
\end{align*}
Applying Cauchy-Schwartz inequality, Young inequality and \eqref{pt}, $I_3$ can be bounded as follows,
\begin{align*}
|I_3|&=2c_0\tau|\sum\limits_{n=1}^{N-1}[(p_t^{n+1}-\partial_\tau p^{n+1},\theta_p^{n+1})+(\partial_\tau\rho_p^{n+1},\theta_p^{n+1})]|\\
&\leq 2c_0\tau\sum\limits_{n=1}^{N-1}||\theta_p^{n+1} ||^2+c_0\tau\sum\limits_{n=1}^{N-1}\|p_t^{n+1}-\partial_\tau p^{n+1}\|^2 +c_0\tau\sum\limits_{n=1}^{N-1}||\partial_\tau\rho_p^{n+1}||^2\\
&\leq 2c_0\tau\sum\limits_{n=1}^{N-1}||\theta_p^{n+1} ||^2+C\tau^2\int_{t_1}^{t_{N}}\|p_{tt}\|^2dt +Ch^{2k_2+2}\int_{t_1}^{t_{N}}\|p_t\|_{k_2+1}^2dt.
\end{align*}
Similarly,
\begin{align*}
|I_4|&=2b_0\tau|\sum\limits_{n=1}^{N-1}[(T_t^{n}-\partial_\tau T^{n},\theta_p^{n+1})+(\partial_\tau\rho_T^{n},\theta_p^{n+1})]|\\
&\leq 2b_0\tau\sum\limits_{n=1}^{N-1}||\theta_p^{n+1} ||^2+b_0\tau\sum\limits_{n=1}^{N-1}\|T_t^{n}-\partial_\tau T^{n}\|^2 +b_0\tau\sum\limits_{n=1}^{N-1}||\partial_\tau\rho_T^{n}||^2\\
&\leq 2b_0\tau\sum\limits_{n=1}^{N-1}||\theta_p^{n+1} ||^2+C\tau^2\int_{t_{0}}^{t_{N-1}}\|T_{tt}\|^2dt +Ch^{2k_3+2}\int_{t_{0}}^{t_{N-1}}\|T_t\|_{k_3+1}^2dt,\\
|I_5|&=2a_0\tau|\sum\limits_{n=1}^{N-1}[(T_t^{n+1}-\partial_\tau T^{n+1},\theta_T^{n+1})+(\partial_\tau\rho_T^{n+1},\theta_T^{n+1})]|\\
&\leq 2a_0\tau\sum\limits_{n=1}^{N-1}||\theta_T^{n+1} ||^2+C\tau^2\int_{t_{1}}^{t_{N}}\|T_{tt}\|^2dt +Ch^{2k_3+2}\int_{t_{1}}^{t_{N}}\|T_t\|_{k_3+1}^2dt,\\
|I_6|&=2b_0\tau|\sum\limits_{n=1}^{N-1}[(p_t^{n+1}-\partial_\tau p^{n+1},\theta_T^{n+1})+(\partial_\tau\rho_p^{n+1},\theta_T^{n+1})]|\\
&\leq 2b_0\tau\sum\limits_{n=1}^{N-1}||\theta_T^{n+1} ||^2+C\tau^2\int_{t_1}^{t_{N}}\|p_{tt}\|^2dt +Ch^{2k_2+2}\int_{t_1}^{t_{N}}\|p_t\|_{k_2+1}^2dt.
\end{align*}
With respect to $I_7$, applying Cauchy-Schwartz inequality, Young's inequality and Lemma \ref{lem2.2}, we derive
\begin{align*}
&\alpha b({\mathbf u}_t^{n+1}-{\mathbf u}_t^n,\theta_p^{n+1})\notag\\
=&\alpha\sum\limits_{K\in \mathcal T_h}\int_K \nabla\cdot({\mathbf u}_t^{n+1}-{\mathbf u}_t^n)\theta_p^{n+1}dK-\alpha\sum\limits_{e\in \Gamma_I}\int_{e}\{\theta_p^{n+1}\}[{\mathbf u}_t^{n+1}-{\mathbf u}_t^n]\cdot{\bm n_e}ds\notag\\
\le& \frac{C_3}{24}||\theta_p^{n+1}||^2+\frac{6\alpha^2}{C_3}||\nabla\cdot({\mathbf u}_t^{n+1}-{\mathbf u}_t^n)||^2+\frac{C_3}{24}\sum\limits_{e\in \Gamma_I}h_e||\{\theta_p^{n+1}\}||_{0,e}^2+\frac{6\alpha^2}{C_3}\sum\limits_{e\in \Gamma_I}h_e^{-1}||[{\mathbf u}_t^{n+1}-{\mathbf u}_t^n]||_{0,e}^2\notag\\
\le& \frac{C_3}{24}|\!|\!|\theta_p^{n+1} |\!|\!|^2+C|\!|\!|{\mathbf u}_t^{n+1}-{\mathbf u}_t^n|\!|\!|_V^2.
\end{align*}
Similarly,
\begin{align*}
\alpha b({\mathbf u}_t^{n}-\partial_\tau {\mathbf u}^n,\theta_p^{n+1})
\le& \frac{C_3}{24}|\!|\!|\theta_p^{n+1} |\!|\!|^2+C|\!|\!|{\mathbf u}_t^{n}-\partial_\tau {\mathbf u}^n|\!|\!|_V^2,\\
\alpha b(\partial_\tau\bm\rho_u^n,\theta_p^{n+1})
\le& \frac{C_3}{24}|\!|\!|\theta_p^{n+1} |\!|\!|^2+C||\partial_\tau\bm\rho_u^n||^2++Ch^2||\partial_\tau\nabla\bm\rho_u^n||^2,
\end{align*}
combined with
\begin{equation*}
|\!|\!|\partial_t{\mathbf u}^{n+1}-\partial_t {\mathbf u}^n|\!|\!|_V^2\leq C\tau\int_{t_{n}}^{t_{n+1}}|\!|\!|\partial_{tt}{\mathbf u}|\!|\!|_V^2dt,
\end{equation*}
\begin{equation*}
|\!|\!|\partial_t{\mathbf u}^{n}-\partial_\tau {\mathbf u}^n|\!|\!|_V^2\leq C\tau\int_{t_{n-1}}^{t_n}|\!|\!|\partial_{tt}{\mathbf u}|\!|\!|_V^2dt,
\end{equation*}
\begin{equation*}
||\partial_\tau\bm\rho_u^n||^2\leq C\frac{h^{2k_1+2}}{\tau}\int_{t_{n-1}}^{t_n}\|\partial_{t}{\mathbf u}\|_{k_1+1}^2dt,
\end{equation*}
we provide
\begin{align*}
|I_7|&=2\alpha\tau|\sum\limits_{n=1}^{N-1}[ b({\mathbf u}_t^{n+1}-{\mathbf u}_t^n,\theta_p^{n+1})+ b({\mathbf u}_t^{n}-\partial_\tau {\mathbf u}^{n},\theta_p^{n+1})+ b(\partial_\tau\bm\rho_u^{n},\theta_p^{n+1})]|\\
&\leq \frac{C_3}{4}\tau\sum\limits_{n=1}^{N-1}|\!|\!|\theta_p^{n+1} |\!|\!|^2+C(\tau^2\int_{t_{0}}^{t_{N}}|\!|\!|{\mathbf u}_{tt}|\!|\!|_V^2dt+h^{2k_1+2}\int_{t_{0}}^{t_{N-1}}\|{\mathbf u}_t\|_{k_1+1}^2dt).
\end{align*}
Similarly,
 \begin{align*}
|I_8|&=2\beta\tau|\sum\limits_{n=1}^{N-1}[ b({\mathbf u}_t^{n+1}-{\mathbf u}_t^n,\theta_T^{n+1})+ b({\mathbf u}_t^n-\partial_\tau {\mathbf u}^n,\theta_T^{n+1})+ b(\partial_\tau\bm\rho_u^n,\theta_T^{n+1})]|\\
&\leq \frac{C_3}{4}\tau\sum\limits_{n=1}^{N-1}|\!|\!|\theta_T^{n+1} |\!|\!|^2+C\left(\tau^2\int_{t_{0}}^{t_{N}}|\!|\!|{\mathbf u}_{tt}|\!|\!|_V^2dt+h^{2k_1+2}\int_{t_{0}}^{t_{N-1}}\|{\mathbf u}_t\|_{k_1+1}^2dt\right).
\end{align*}
Applying Cauchy-Schwartz inequality, Young's inequality and Lemma \ref{lem0},
 \begin{align*}
|I_9|=2\tau|\sum\limits_{n=1}^{N-1}(\mathcal M({\mathbf K}\nabla p_h^{n})\cdot\nabla\theta_T^{n+1},\theta_T^{n+1})|\leq \frac{C_3}{8}\tau\sum\limits_{n=1}^{N-1}|\!|\!|\theta_T^{n+1}|\!|\!|^2+C\tau\sum\limits_{n=1}^{N-1}\|\theta_T^{n+1}\|^2.
\end{align*}
Since we suppose that the constant $M$ for the cut-off operator is sufficiently large, by virtue of \cite{MR2116915}, we rewrite $I_{10}$ as follows,
\begin{align*}
&(({\mathbf K}\nabla p^{n+1}-\mathcal M({\mathbf K}\nabla p_h^{n}))\cdot\nabla R_TT^{n+1},\theta_T^{n+1})\\
=&({\mathbf K}\nabla\rho_p^{n+1}\cdot\nabla R_TT^{n+1},\theta_T^{n+1})+({\mathbf K}\nabla(R_pp^{n+1}-R_pp^{n})\cdot\nabla R_TT^{n+1},\theta_T^{n+1})\\
&+((\mathcal M({\mathbf K}\nabla R_pp^{n})-\mathcal M({\mathbf K}\nabla p_h^{n}))\cdot\nabla R_TT^{n+1},\theta_T^{n+1})\\
:=&H_1+H_2+H_3.
\end{align*}
For $H_1$, utilizing integrating by parts and the boundary condition, we get
\begin{align*}
|H_1|\le& k_M|(\nabla T^{n+1}\nabla\rho_p^{n+1},\theta_T^{n+1})-(\nabla \rho_T^{n+1}\cdot\nabla\rho_p^{n+1},\theta_T^{n+1})|\notag\\
=&k_M|\langle\rho_p^{n+1},\nabla T^{n+1}\cdot\theta_T^{n+1}\bm n\rangle_{\partial\Omega}-(\rho_p^{n+1},\nabla T^{n+1}\cdot\nabla\theta_T^{n+1})\notag\\
&-(\rho_p^{n+1},\nabla\cdot(\nabla T^{n+1})\cdot\theta_T^{n+1})-(\nabla \rho_T^{n+1}\cdot\nabla\rho_p^{n+1},\theta_T^{n+1})|\notag\\
=&k_M|(\rho_p^{n+1},\nabla T^{n+1}\cdot\nabla\theta_T^{n+1})+(\rho_p^{n+1},\nabla\cdot(\nabla T^{n+1})\cdot\theta_T^{n+1})\notag\\
&+(\nabla \rho_T^{n+1}\cdot\nabla\rho_p^{n+1},\theta_T^{n+1})|.
\end{align*}
Applying Cauchy-Schwartz inequality, Young's inequality and \eqref{pt0}, we derive
\begin{align*}
(\rho_p^{n+1},\nabla T^{n+1}\cdot\nabla\theta_T^{n+1})\le& \frac{C_3}{8}||\nabla\theta_T^{n+1}||^2+C||\rho_p^{n+1}||^2,\\
(\rho_p^{n+1},\nabla\cdot(\nabla T^{n+1})\cdot\theta_T^{n+1})\le& C||\theta_T^{n+1}||^2+C||\rho_p^{n+1}||^2,\\
(\nabla \rho_T^{n+1}\cdot\nabla\rho_p^{n+1},\theta_T^{n+1})\le& C||\nabla(R_hp^{n+1}-R_hp^{n}||^2+C||\theta_T^{n+1}||^2\\
\le& C\tau\int_{t_{n}}^{t_{n+1}}||\nabla\partial_tR_hp||^2dt+C||\theta_T^{n+1}||^2.
\end{align*}
By virtue of the $L^{\infty}(\Omega)$-error estimate in \cite{MR717695}, we have
\begin{equation*}
||\nabla\rho_T^{n+1}||_{0,\infty}\le C||T||_{k_3+1,\infty}h^{k_3}.
\end{equation*}
Therefore,
\begin{align*}
(\nabla \rho_T^{n+1}\cdot\nabla\rho_p^{n+1},\theta_T^{n+1})
&\le\frac{1}{2}||\theta_T^{n+1}||^2+\frac{1}{2}||\nabla \rho_T^{n+1}||_{0,\infty}^2||\nabla\rho_p^{n+1}||^2\\
&\le C||\theta_T^{n+1}||^2+Ch^{2k_2}h^{2k_3}\\
&\le C||\theta_T^{n+1}||^2+C(h^{2k_2+2}+h^{2k_3+2}).
\end{align*}
Then, we have
\begin{align*}
|H_1|\le C(h^{2k_2+2}+h^{2k_3+2})+\frac{C_3}{8}||\nabla\theta_T^{n+1}||^2+C||\theta_T^{n+1}||^2.
\end{align*}
By \eqref{pt0}, we obtain
\begin{align*}
|H_2|\le C||\nabla (R_pp^{n+1}-R_pp^{n})||^2+C||\theta_T^{n+1}||^2\le C\tau\int_{t_{n}}^{t_{n+1}}||\nabla R_pp_t||^2dt+C||\theta_T^{n+1}||^2.
\end{align*}
Noting that $|\mathcal M({\mathbf K}\nabla p^{n+1})-\mathcal M({\mathbf K}\nabla p_h^{n})|\le|{\mathbf K}\nabla p^{n+1}-{\mathbf K}\nabla p_h^{n}|$, hence
\begin{align*}
|H_3|\le||\nabla R_TT^{n+1}||_{0,\infty} ||{\mathbf K}\nabla R_p p^{n}-{\mathbf K}\nabla p_h^{n}||||\theta_T^{n+1}||\le \frac{C_3}{4}||\nabla\theta_p^{n}||^2+C||\theta_T^{n+1}||^2.
\end{align*}
By combining $H_1$-$H_3$, we obtain
\begin{align*}
&(({\mathbf K}\nabla p^{n+1}-\mathcal M({\mathbf K}\nabla p_h^{n}))\cdot\nabla R_TT^{n+1},\theta_T^{n+1})\\
\leq & C(h^{2k_2+2}+h^{2k_3+2})+C\tau\int_{t_{n}}^{t_{n+1}}||\nabla p_t||_2^2dt+\frac{C_3}{4}\|\nabla\theta_p^{n}\|^2+\frac{C_3}{8}\|\nabla\theta_T^{n+1}\|^2+C||\theta_T^{n+1}||^2.
\end{align*}
Thus,
\begin{align*}
|I_{10}|=&2\tau|\sum\limits_{n=1}^{N-1}(({\mathbf K}\nabla p^{n+1}-\mathcal M({\mathbf K}\nabla p_h^{n}))\cdot\nabla R_TT^{n+1},\theta_T^{n+1})|\\
\leq &CT(h^{2k_2+2}+h^{2k_3+2})+C\tau^2\int_{t_{1}}^{t_{N}}||\nabla p_t||_2^2dt\\
&+\frac{C_3}{4}\tau\sum\limits_{n=1}^{N-1}\|\nabla\theta_p^{n}\|^2+\frac{C_3}{8}\tau\sum\limits_{n=1}^{N-1}\|\nabla\theta_T^{n+1}\|^2+C\tau\sum\limits_{n=1}^{N-1}||\theta_T^{n+1}||^2.
\end{align*}
Similarly to $I_7$, estimated $I_{11}$ and $I_{12}$,
\begin{align*}
&\alpha(b(\partial_\tau\bm\theta_u^{n+1},\theta_p^{n+1})-b(\partial_\tau\bm\theta_u^n,\theta_p^{n+1}))\\
=&\alpha\sum\limits_{K\in \mathcal T_h}\int_K (\partial_\tau\bm\theta_u^{n+1}-\partial_\tau\bm\theta_u^n)\cdot\nabla\theta_p^{n+1}dK+\alpha\sum\limits_{e\in \Gamma_I}\int_{e}\{\partial_\tau\bm\theta_u^{n+1}-\partial_\tau\bm\theta_u^n\}\cdot{\bm n_e}[\theta_p^{n+1}]ds\\
\le&\alpha||\partial_\tau\bm\theta_u^{n+1}-\partial_\tau\bm\theta_u^n||||\nabla\theta_p^{n+1}||
+\alpha\sum\limits_{e\in \Gamma_I}||\{\partial_\tau\bm\theta_u^{n+1}-\partial_\tau\bm\theta_u^n\}||_{0,e}||[\theta_p^{n+1}]||_{0,e}\\
\le& \frac{C_3}{2}|\!|\!|\theta_p^{n+1}|\!|\!|^2+\frac{\alpha^2}{2C_3}(1+C)||\partial_\tau\bm\theta_u^{n+1}-\partial_\tau\bm\theta_u^n||^2,
\end{align*}
we obtain
\begin{align*}
&|I_{11}|+|I_{12}|\\
=&2\alpha\tau|\sum\limits_{n=1}^{N-1}(b(\partial_\tau\bm\theta_u^{n+1},\theta_p^{n+1})-b(\partial_\tau\bm\theta_u^n,\theta_p^{n+1}))|+2\beta\tau|\sum\limits_{n=1}^{N-1}(b(\partial_\tau\bm\theta_u^{n+1},\theta_T^{n+1})-b(\partial_\tau\bm\theta_u^n,\theta_T^{n+1}))|\\
\le& C_3\tau\sum\limits_{n=1}^{N-1}|\!|\!|\theta_p^{n+1}|\!|\!|^2+ C_3\tau\sum\limits_{n=1}^{N-1}|\!|\!|\theta_T^{n+1}|\!|\!|^2+\frac{\alpha^2+\beta^2}{C_3}(1+C)\tau\sum\limits_{n=1}^{N-1}||\partial_\tau\bm\theta_u^{n+1}-\partial_\tau\bm\theta_u^n||^2.
\end{align*}
From \eqref{f}, we get
\begin{align*}
I_{13}=2\alpha\tau\sum\limits_{n=1}^{N-1}  b(\partial_\tau\bm\theta_u^{n+1},\rho_p^{n+1})=2\alpha b(\bm\theta_u^{N},\rho_p^{N}) - 2\alpha b(\bm\theta_u^{1},\rho_p^{1}) -2\alpha\sum\limits_{n=1}^{N-1}  b(\bm\theta_u^{n},\rho_p^{n+1}-\rho_p^{n}).
\end{align*}
Applying Cauchy-Schwartz inequality and Young's inequality yields that
\begin{align*}
\alpha b(\bm\theta_u^{N},\rho_p^{N})
=&\alpha\sum\limits_{K\in \mathcal T_h}\int_K \rho_p^{N}\nabla\cdot \bm\theta_u^{N}dK-\alpha \sum\limits_{e\in \Gamma}\int_{e}[\bm\theta_u^n]\cdot{\bm n_e}\{\rho_p^{N}\}ds\\
\le& \alpha\delta|\!|\!|\bm\theta_u^{N}|\!|\!|_V^2+C||\rho_p^{N}||^2+Ch^2||\nabla\rho_p^{N}||^2\\
\le& \frac{C_1}{6}|\!|\!|\bm\theta_u^{N}|\!|\!|_V^2+Ch^{2k_2+2},\\
\alpha b(\bm\theta_u^{1},\rho_p^{1})
\le& C|\!|\!|\bm\theta_u^{1}|\!|\!|_V^2+Ch^{2k_2+2},\\
\alpha \sum\limits_{n=1}^{N-1}  b(\bm\theta_u^{n},\rho_p^{n+1}-\rho_p^{n})
\le&\alpha C_5\tau\sum\limits_{n=1}^{N-1} |\!|\!|\bm\theta_u^{n}|\!|\!|_V^2+\frac{\alpha C_5}{\tau}\sum\limits_{n=1}^{N-1}||\rho_p^{n+1}-\rho_p^{n}||^2\\
\le& \alpha C\tau\sum\limits_{n=1}^{N-1}|\!|\!|\bm\theta_u^{n}|\!|\!|_V^2+C\int_{t_1}^{t_N}||\frac{\partial\rho_p}{\partial t}||^2dt\\
\le& C\tau\sum\limits_{n=1}^{N-1}|\!|\!|\bm\theta_u^{n}|\!|\!|_V^2+Ch^{2k_2+2}.
\end{align*}
Hence,
\begin{align*}
|I_{13}|
\leq \frac{C_1}{6}|\!|\!|\bm\theta_u^{N}|\!|\!|_V^2+C|\!|\!|\bm\theta_u^{1}|\!|\!|_V^2+ C\tau\sum\limits_{n=1}^{N-1} |\!|\!|\bm\theta_u^{n}|\!|\!|_V^2+Ch^{2k_2+2}.
\end{align*}
Similarly,
\begin{align*}
|I_{14}|
\leq \frac{C_1}{6}|\!|\!|\bm\theta_u^{N}|\!|\!|_V^2+C|\!|\!|\bm\theta_u^{1}|\!|\!|_V^2+ C\tau\sum\limits_{n=1}^{N-1} |\!|\!|\bm\theta_u^{n}|\!|\!|_V^2+Ch^{2k_3+2}.
\end{align*}
In addition, applying \eqref{f}, we have
\begin{align*}
&I_{15}=\sum\limits_{n=1}^{N-1}\gamma(\partial_\tau R_u{\mathbf u}^{n+1}-\partial_\tau R_u{\mathbf u}^{n},\bm\theta_u^{n+1}-\bm\theta_u^{n})\\
=&\gamma(\frac{R_u{\mathbf u}^{2}-R_u{\mathbf u}^{1}}{\tau}-\frac{R_u{\mathbf u}^{1}-R_u{\mathbf u}^{0}}{\tau_0},\bm\theta_u^{2}-\bm\theta_u^{1})+\gamma(\frac{R_u{\mathbf u}^{3}-R_u{\mathbf u}^{2}}{\tau}-\frac{R_u{\mathbf u}^{2}-R_u{\mathbf u}^{1}}{\tau},\bm\theta_u^{3}-\bm\theta_u^{2})\\
&+\sum\limits_{n=3}^{N-1}\gamma(\frac{R_u{\mathbf u}^{n+1}-R_u{\mathbf u}^{n}}{\tau}-\frac{R_u{\mathbf u}^{n}-R_u{\mathbf u}^{n-1}}{\tau},\bm\theta_u^{n+1}-\bm\theta_u^{n})\\
=&\gamma(\frac{R_u{\mathbf u}^{2}-R_u{\mathbf u}^{1}}{\tau}-\frac{R_u{\mathbf u}^{1}-R_u{\mathbf u}^{0}}{\tau_0},\bm\theta_u^{2}-\bm\theta_u^{1})+\gamma(\frac{R_u{\mathbf u}^{3}-R_u{\mathbf u}^{2}}{\tau}-\frac{R_u{\mathbf u}^{2}-R_u{\mathbf u}^{1}}{\tau},\bm\theta_u^{3}-\bm\theta_u^{2})\\
&+\gamma(\frac{R_u{\mathbf u}^{N}-R_u{\mathbf u}^{N-1}}{\tau}-\frac{R_u{\mathbf u}^{N-1}-R_u{\mathbf u}^{N-2}}{\tau},\bm\theta_u^{N})-\gamma(\frac{R_u{\mathbf u}^{3}-R_u{\mathbf u}^{2}}{\tau}-\frac{R_u{\mathbf u}^{2}-R_u{\mathbf u}^{1}}{\tau},\bm\theta_u^{3})\\
&-\sum\limits_{n=3}^{N-1}\gamma(\frac{R_u{\mathbf u}^{n+1}-2R_u{\mathbf u}^{n}+R_u{\mathbf u}^{n-1}}{\tau}-\frac{R_u{\mathbf u}^{n}-2R_u{\mathbf u}^{n-1}+R_u{\mathbf u}^{n-2}}{\tau},\bm\theta_u^{n})\\
:=&J_1+J_2+J_3+J_4+J_5.
\end{align*}
Applying Cauchy-Schwartz inequality and Young's inequality, we get
\begin{align*}
J_1=&\gamma(\frac{R_u{\mathbf u}^{2}-R_u{\mathbf u}^{1}}{\tau}-\frac{R_u{\mathbf u}^{1}-R_u{\mathbf u}^{0}}{\tau_0},\bm\theta_u^{2}-\bm\theta_u^{1})\\
=&\gamma(\frac{1}{\tau}\int_{t_{1}}^{t_{2}}(t_2-t)R_u{\mathbf u}_{tt}dt+\frac{1}{\tau}\int_{t_{0}}^{t_{1}}(t-t_0)R_u{\mathbf u}_{tt}dt,\bm\theta_u^{2}-\bm\theta_u^{1})\\
\le&\gamma||\int_{t_{0}}^{t_{2}}R_u{\mathbf u}_{tt}dt||||\bm\theta_u^{2}-\bm\theta_u^{1}||\le\frac{\gamma}{4\tau}||\bm\theta_u^{2}-\bm\theta_u^{1}||^2+C\tau^2\int_{t_{0}}^{t_{2}}||R_u{\mathbf u}_{tt}||^2dt,\\
J_2=&\gamma(\frac{R_u{\mathbf u}^{3}-R_u{\mathbf u}^{2}}{\tau}-\frac{R_u{\mathbf u}^{2}-R_u{\mathbf u}^{1}}{\tau},\bm\theta_u^{3}-\bm\theta_u^{2})\le\frac{\gamma}{4\tau}||\bm\theta_u^{3}-\bm\theta_u^{2}||^2+C\tau^2\int_{t_{1}}^{t_{3}}||R_u{\mathbf u}_{tt}||^2dt,\\
J_3=&\gamma(\frac{R_u{\mathbf u}^{N}-R_u{\mathbf u}^{N-1}}{\tau}-\frac{R_u{\mathbf u}^{N-1}-R_u{\mathbf u}^{N-2}}{\tau},\bm\theta_u^{N})\\
=&\gamma(\frac{1}{\tau}\int_{t_{N-1}}^{t_{N}}(t_N-t)R_u{\mathbf u}_{tt}dt+\frac{1}{\tau}\int_{t_{N-2}}^{t_{N-1}}(t-t_{N-2})R_u{\mathbf u}_{tt}dt,\bm\theta_u^{N})\\
\le&\gamma||\int_{t_{N-2}}^{t_{N}}R_u{\mathbf u}_{tt}dt||||\bm\theta_u^{N}||\le C\tau(\max\limits_{0\le t\le t_f}||R_u{\mathbf u}_{tt}||)|\!|\!|\bm\theta_u^{N}|\!|\!|_V\\
\le&\frac{C_1}{6}|\!|\!|\bm\theta_u^{N}|\!|\!|_V^2+C\tau^2\max\limits_{0\le t\le t_f}||R_u{\mathbf u}_{tt}||^2,\\
J_4=&-\gamma(\frac{R_u{\mathbf u}^{3}-R_u{\mathbf u}^{2}}{\tau}-\frac{R_u{\mathbf u}^{2}-R_u{\mathbf u}^{1}}{\tau},\bm\theta_u^{3})\le\frac{C_1}{4}|\!|\!|\bm\theta_u^{3}|\!|\!|_V^2+C\tau^2\max\limits_{0\le t\le t_f}||R_u{\mathbf u}_{tt}||^2.
\end{align*}
From the Taylor expansion, it follows that
\begin{align*}
\frac{R_u{\mathbf u}^{n+1}-2R_u{\mathbf u}^{n}+R_u{\mathbf u}^{n-1}}{\tau}
&=\tau R_u{\mathbf u}_{tt}^{n}+\frac{1}{2\tau}\int_{t_n}^{t_{n+1}}(t-t_{n+1})^2R_u{\mathbf u}_{ttt}dt-\frac{1}{2\tau}\int_{t_{n-1}}^{t_{n}}(t-t_{n-1})^2R_u{\mathbf u}_{ttt}dt,\\
\frac{R_u{\mathbf u}^{n}-2R_u{\mathbf u}^{n-1}+R_u{\mathbf u}^{n-2}}{\tau}
&=\tau R_u{\mathbf u}_{tt}^{n-1}+\frac{1}{2\tau}\int_{t_{n-1}}^{t_{n}}(t-t_{n})^2R_u{\mathbf u}_{ttt}dt-\frac{1}{2\tau}\int_{t_{n-2}}^{t_{n-1}}(t-t_{n-2})^2R_u{\mathbf u}_{ttt}dt.
\end{align*}
Therefore, we get
\begin{align*}
&||\frac{R_u{\mathbf u}^{n+1}-2R_u{\mathbf u}^{n}+R_u{\mathbf u}^{n-1}}{\tau}-\frac{R_u{\mathbf u}^{n}-2R_u{\mathbf u}^{n-1}+R_u{\mathbf u}^{n-2}}{\tau}||\\
=&||\tau\int_{t_{n-1}}^{t_n}R_u{\mathbf u}_{ttt}dt+\frac{1}{2\tau}\int_{t_n}^{t_{n+1}}(t-t_{n+1})^2R_u{\mathbf u}_{ttt}dt-\frac{1}{2\tau}\int_{t_{n-1}}^{t_{n}}(t-t_{n-1})^2R_u{\mathbf u}_{ttt}dt\\
&-\frac{1}{2\tau}\int_{t_{n-1}}^{t_{n}}(t-t_{n})^2R_u{\mathbf u}_{ttt}dt+\frac{1}{2\tau}\int_{t_{n-2}}^{t_{n-1}}(t-t_{n-2})^2R_u{\mathbf u}_{ttt}dt||\\
\le&2\tau\int_{t_{n-2}}^{t_{n+1}}||R_u{\mathbf u}_{ttt}||dt.
\end{align*}
Applying Cauchy-Schwartz inequality and Young's inequality, we get
\begin{align*}
J_5=&-\sum\limits_{n=3}^{N-1}\gamma(\frac{R_u{\mathbf u}^{n+1}-2R_u{\mathbf u}^{n}+R_u{\mathbf u}^{n-1}}{\tau}-\frac{R_u{\mathbf u}^{n}-2R_u{\mathbf u}^{n-1}+R_u{\mathbf u}^{n-2}}{\tau},\bm\theta_u^{n})\\
\le&\gamma\sum\limits_{n=3}^{N-1}||\frac{R_u{\mathbf u}^{n+1}-2R_u{\mathbf u}^{n}+R_u{\mathbf u}^{n-1}}{\tau}-\frac{R_u{\mathbf u}^{n}-2R_u{\mathbf u}^{n-1}+R_u{\mathbf u}^{n-2}}{\tau}||||\bm\theta_u^{n}||\\
\le&\gamma C\tau\sum\limits_{n=3}^{N-1}|\!|\!|\bm\theta_u^{n}|\!|\!|_V^2+C\tau^2\int_{t_1}^{t_N}||R_u{\mathbf u}_{ttt}||^2dt.
\end{align*}
Combining the estimates from $J_1$-$J_5$, we get
\begin{align*}
|I_{15}|\le&\frac{C_1}{6}|\!|\!|\bm\theta_u^{N}|\!|\!|_V^2+\frac{C_1}{4}|\!|\!|\bm\theta_u^{3}|\!|\!|_V^2+\frac{\gamma}{4\tau}||\bm\theta_u^{2}-\bm\theta_u^{1}||^2+\frac{\gamma}{4\tau}||\bm\theta_u^{3}-\bm\theta_u^{2}||^2\\
&+\gamma C\tau\sum\limits_{n=3}^{N-1}|\!|\!|\bm\theta_u^{n}|\!|\!|_V^2+C\tau^2(\max\limits_{0\le t\le t_f}||R_u{\mathbf u}_{tt}||^2+\int_{t_1}^{t_N}||R_u{\mathbf u}_{ttt}||^2dt).
\end{align*}
Synthesizing $I_1$-$I_{15}$, we obtain
\begin{align}
&(c_0-2b_0)||\theta_p^{N}||^2+(c_0-2b_0)\tau^2\sum\limits_{n=1}^{N-1}||\partial_\tau\theta_p^{n+1}||^2+(a_0-b_0)||\theta_T^{N}||^2\notag\\
&+(a_0-2b_0)\tau^2\sum\limits_{n=1}^{N-1}||\partial_\tau\theta_T^{n+1}||^2+\frac{C_1}{2}|\!|\!|\bm\theta_u^{N}|\!|\!|_V^2+C_1\tau^2\sum\limits_{n=1}^{N-1}|\!|\!|\partial_\tau{\bm\theta}_u^{n+1}|\!|\!|_V^2+ \frac{C_3}{2}\tau\sum\limits_{n=1}^{N-1}|\!|\!|\theta_p^{n+1} |\!|\!|^2\notag\\
&+ \frac{C_3}{2}\tau\sum\limits_{n=1}^{N-1}|\!|\!|\theta_T^{n+1} |\!|\!|^2+\gamma\tau||{\partial_\tau\bm\theta}_u^{N}||^2+\left(\gamma-\frac{\alpha^2+\beta^2}{C_3}(1+C)\right)\tau\sum\limits_{n=1}^{N-1}||{\partial_\tau\bm\theta}_u^{n+1}-{\partial_\tau\bm\theta}_u^{n}||^2\notag\\
\le& C\tau\sum\limits_{n=1}^{N-1}||\theta_p^{n+1} ||^2+C\tau \sum\limits_{n=1}^{N-1}||\theta_T^{n+1}||^2+C\tau\sum\limits_{n=1}^{N-1} |\!|\!|\bm\theta_u^{n}|\!|\!|_V^2+ C(\tau^2+h^{2k_1+2}+h^{2k_2+2}+h^{2k_3+2})\notag\\
&+(c_0+2b_0)||\theta_p^{1}||^2+(a_0+b_0)||\theta_T^{1}||^2+ C|\!|\!|\bm\theta_u^{1}|\!|\!|_V^2 +a({\bm\theta}_u^{1},{\bm\theta}_u^{1})+\gamma\tau||{\partial_\tau\bm\theta}_u^{1}||^2\notag\\
&+\frac{C_1}{4}|\!|\!|\bm\theta_u^{3}|\!|\!|_V^2+\frac{\gamma}{4\tau}||\bm\theta_u^{2}-\bm\theta_u^{1}||^2+\frac{\gamma}{4\tau}||\bm\theta_u^{3}-\bm\theta_u^{2}||^2.\label{z1}
\end{align}
In \eqref{z1}, choosing $N=2$, we get
\begin{align*}
&\frac{\gamma}{4\tau}||\bm\theta_u^{2}-\bm\theta_u^{1}||^2 \\
\le&C(\tau^2+h^{2k_1+2}+h^{2k_2+2}+h^{2k_3+2})\notag\\
&+(c_0+b_0)||\theta_p^{1}||^2+(a_0+b_0)||\theta_T^{1}||^2+ a({\bm\theta}_u^{1},{\bm\theta}_u^{1})+\gamma\tau||{\partial_\tau\bm\theta}_u^{1}||^2.
\end{align*}
In \eqref{z1}, choosing $N=3$ and utilizing Gronwall inequality, we get
\begin{align*}
&\frac{C}{4}|\!|\!|\bm\theta_u^{3}|\!|\!|_V^2+\frac{\gamma}{4\tau}||\bm\theta_u^{3}-\bm\theta_u^{2}||^2\\
\le&C(\tau^2+h^{2k_1+2}+h^{2k_2+2}+h^{2k_3+2})\notag\\
&+(c_0+b_0)||\theta_p^{1}||^2+(a_0+b_0)||\theta_T^{1}||^2+ a({\bm\theta}_u^{1},{\bm\theta}_u^{1})+\gamma\tau||{\partial_\tau\bm\theta}_u^{1}||^2.
\end{align*}
Because of the triangle inequality, Lemma \ref{lem3.7}, Gronwall inequality and noticing that $\bm\theta_u^0= \bm 0,$ choosing $\gamma $ large enough, when $c_0-2b_0>0$, $a_0-2b_0\ge0$ and $\gamma-\frac{\alpha^2+\beta^2}{C_3}(1+C)\ge0,$ we get
\begin{align*}
&|\!|\!|\mathbf u(t_n)-\mathbf u_h^n|\!|\!|_V^2+||p(t_{n})-p_h^n||^2+||T(t_{n})-T_h^n||^2\\
\le &C(\tau^2+h^{2k_1}+h^{2k_2+2}+h^{2k_3+2})+C||p_h^{1}-R_pp^1||^2\\
&+C||T_h^{1}-R_TT^1||^2+C|\!|\!|\mathbf u_h^{1}-R_u\mathbf u^1|\!|\!|_V^2+\gamma\tau||\frac{\mathbf u_h^{1}-R_u\mathbf u^1}{\tau_0}||^2,
\end{align*}
which finishes the proof.
\end{proof}
\begin{remark}\label{remark2}
From Theorem \ref{thm1}, it can be seen that the results depend on $\mathbf u_h^1, p_h^1$ and $T_h^1$. If the approximate error at the initial time step is optimal, the proposed SDG method (\ref{s1})-(\ref{s3}) can achieve the optimal convergence speed.
\end{remark}

\begin{lemma}\label{lem001}\cite{chen}
Let the initial values satisfy (A1)-(A6) and $({\mathbf u}_h^1,p_h^1,T_h^1)\in {\mathbf V}_h\times Q_h\times S_h$ be the solution of (\ref{b71})-(\ref{b91}). Then, we have the following error estimates
\begin{align*}
|\!|\!|\mathbf u_h^{1}-R_u\mathbf u^1|\!|\!|_V^2+||p_h^{1}-R_pp^1||^2+||T_h^{1}-R_TT^1||^2&\leq C(\tau^2+h^{2k_1}+h^{2k_2+2}+h^{2k_3+2}),\\
\tau|\!|\!|\frac{\mathbf u_h^{1}-R_u\mathbf u^1}{\tau}|\!|\!|_V^2&\leq C(\tau^2+h^{2k_1}+h^{2k_2+2}+h^{2k_3+2}).
\end{align*}
\end{lemma}
\begin{cor}\label{cor}
The numerical solution of (\ref{s1})-(\ref{s3}) satisfies the following error estimate,
\begin{align*}
|\!|\!|\mathbf u^n-\mathbf u_h^n|\!|\!|_V+||p^{n}-p_h^n||+||T^{n}-T_h^n||\le C(\tau+h^{k_1}+h^{k_2+1}+h^{k_3+1}).
\end{align*}
\end{cor}
\begin{proof}
Lemma \ref{lem001} and Theorem \ref{thm1} directly lead to the conclusion.
\end{proof}

\section{Numerical experiments}\label{NE}
In this section, we present numerical tests under the various coupling strengths to validate our theoretical results.  Note that the calculations performed here are powered by the State Key Laboratory of Scientific and Engineering Computing, Chinese Academy of Sciences, high performance computers.  Each computer node has two $18$-core $2.3$ GHz Intel Xeon Gold $6140$ processors and $192$ GB of RAM.  The package PETSc \cite{petsc-web-page,petsc-user-ref,petsc-efficient} solves the linear equations created by the coefficient matrix from \eqref{e1}-\eqref{e3}.

With the homogeneous Dirichlet boundary conditions \eqref{e4} for $\mathbf u$, $p$ and $T$, we examine the system \eqref{e1}-\eqref{e3} on a two dimensional domain $\Omega=[0,1]^2$ in our example. The parameters we take are $\gamma=5$, $\lambda=1$, $\mu=1$, $\mathbf K=\begin{pmatrix} 1 & 0 \\ 0 & 1 \end{pmatrix}$ and $\bm\Theta=\begin{pmatrix} 1 & 0 \\ 0 & 1 \end{pmatrix}$. We set $\alpha$, $\beta$, and $b_0$ as shown in Table \ref{table1} with $a_0=c_0=3b_0$, which shows five distinct parameter choices, PA1-PA5, to weaken or strengthen the coupling between the equations in order to test our algorithm under the various coupling strengths.  Every term on the right-hand side has been chosen based on the analytical solution
$$
\displaystyle\mathbf u=\binom{\exp^{-t}\sin(\pi x)\sin(\pi y)}{\exp^{-t}\sin(\pi x)\sin(\pi y)},
$$
$$
p=t\sin(\pi x)\sin(\pi y),
$$
$$
T=\exp(-t)\sin(\pi x)\sin(\pi y).
$$
We select $k_1=k_2=k_3=1$ and employ a uniform triangular mesh $\mathcal{T}_h$ for the discontinuous finite element spaces, the final time $t_f=1$ and the two penalty parameters $\sigma_1=\sigma_2=1\textrm{e}$+$6$. Regarding the calculation of initial values, we adopt \textbf{Option 2}.

\begin{table}[!htbp]
\caption{Selection of parameters for the strong or weak coupling among equations.}\label{table1}
  \begin{center}
  \begin{tabular}{cccccc}
   \hline
        &          PA1&       PA2&       PA3&          PA4&           PA5 \\ \hline
$\alpha$&          1.0&       0.1&       0.1&          1.0&           0.1 \\
$\beta$ &          1.0&       0.1&       1.0&          0.1&           0.1 \\
$b_0$   &          1.0&       1.0&       0.1&          0.1&           0.1 \\
\hline
  \end{tabular}
  \end{center}
\end{table}

In order to observe the temporal convergence order, we fix the spatial step size $h=\frac{1} {128} $, gradually subdividing the time step size. Table \ref {table16} shows the time error convergence order of displacement, pressure and temperature under parameter PA1, which
 confirms the first order convergence.

\begin{table}[htpb]
         \centering
	\caption{ Convergence order of temporal error under PA1 parameter}
          \label{table16}
            \begin{center}
  \begin{tabular}{ccccccccc}
   \hline
  $\tau$&   $||{\mathbf u}(t_n)-{\mathbf u}_h^n||$&    R&   $||p(t_n)-p_h^n||$&         R&    $||T(t_n)-T_h^n||$&          R \\   \hline
$\frac{1}{4}$&             1.39629e-02&             -&         9.69476e-03&               -&               2.04177e-02&   -\\
$\frac{1}{8}$&             3.81250e-03&       1.8728 &         3.02490e-03&         1.6803&                1.06271e-02 &       0.9421\\
$\frac{1}{16}$&            2.32912e-03&             0.7109&    1.71721e-03&               0.8168&          5.55618e-03&                  0.9356\\
$\frac{1}{32}$&           1.05260e-03&             1.1458&    8.25810e-04&               1.0562&          2.81893e-03&                  0.9789\\
$\frac{1}{64}$&           5.39123e-04&             0.9653&    4.25382e-04&               0.9571&          1.43728e-03&     0.9718\\
   \hline
          \end{tabular}
            \end{center}
        \end{table}

Table \ref{table11}-\ref{table15} show the spatial convergence orders for error with time step $\tau=h^2$. The theoretical conclusions of Corollary \ref{cor} are validated by these tables, which show that our proposed method can reach the optimal convergence orders for error within each of the five different coupling strengths.
\begin{table}[htpb]
         \centering
	\caption{ Convergence order of spatial error under PA1 parameter}
          \label{table11}
            \begin{center}
  \begin{tabular}{ccccccccc}
   \hline
  $h$&   $|\!|\!|{\mathbf u}(t_n)-{\mathbf u}_h^n|\!|\!|_V$&    R&   $||p(t_n)-p_h^n||$&         R&    $||T(t_n)-T_h^n||$&          R \\    \hline
$\frac{1}{4}$&             5.11013e-01&             -&         7.35984e-02&               -&               4.62957e-02&   -\\
$\frac{1}{8}$&             2.60434e-01&       0.9724 &         1.93593e-02&               1.9266&          1.32480e-02 &       1.8051\\
$\frac{1}{16}$&            1.30818e-01&             0.9934&    4.90354e-03&               1.9811&          3.42308e-03&                  1.9524\\
$\frac{1}{32}$&           6.54833e-02&             0.9984&    1.22995e-03&               1.9953&          8.62793e-04&                  1.9882\\
\hline
  \end{tabular}
  \end{center}
\end{table}

\begin{table}[htpb]
         \centering
	\caption{ Convergence order of spatial error under PA2 parameter}
          \label{table12}
            \begin{center}
  \begin{tabular}{ccccccccc}
   \hline
  $h$&   $|\!|\!|{\mathbf u}(t_n)-{\mathbf u}_h^n|\!|\!|_V$&    R&   $||p(t_n)-p_h^n||$&         R&    $||T(t_n)-T_h^n||$&          R \\    \hline
$\frac{1}{4}$&             5.10316e-01&             -&         7.37022e-02&               -&               4.64834e-02&   -\\
$\frac{1}{8}$&             2.60308e-01&       0.9712 &  1.94079e-02&               1.9251&          1.33347e-02 &       1.8015\\
$\frac{1}{16}$&            1.30801e-01&      0.9928&    4.91748e-03&               1.9807&          3.447960e-03&        1.9514\\
$\frac{1}{32}$&            6.54811e-02&       0.9982&   1.23350e-03&               1.9952&          8.692260e-04&         1.9879\\
\hline
  \end{tabular}
  \end{center}
\end{table}

\begin{table}[htpb]
         \centering
	\caption{ Convergence order of spatial error under PA3 parameter}
          \label{table13}
            \begin{center}
  \begin{tabular}{ccccccccc}
   \hline
$h$&   $|\!|\!|{\mathbf u}(t_n)-{\mathbf u}_h^n|\!|\!|_V$&    R&   $||p(t_n)-p_h^n||$&         R&           $||T(t_n)-T_h^n||$&          R \\    \hline
$\frac{1}{4}$&            5.10423e-01&             -&   7.85425e-02&               -&               4.19485e-02&   -\\
$\frac{1}{8}$&            2.60328e-01&       0.9714 &    2.09593e-02&               1.9059&         1.18891e-02 &       1.8190\\
$\frac{1}{16}$&           1.30803e-01&      0.9929&    5.33105e-03&               1.9751&           3.06263e-03&        1.9568\\
$\frac{1}{32}$&           6.54814e-02&       0.9982&    1.33865e-03&               1.9936&          7.71287e-04&         1.9894\\
\hline
  \end{tabular}
  \end{center}
\end{table}

\begin{table}[htpb]
         \centering
	\caption{ Convergence order of spatial error under PA4 parameter}
          \label{table14}
            \begin{center}
  \begin{tabular}{ccccccccc}
   \hline
$h$&   $|\!|\!|{\mathbf u}(t_n)-{\mathbf u}_h^n|\!|\!|_V$&    R&   $||p(t_n)-p_h^n||$&         R&           $||T(t_n)-T_h^n||$&          R \\    \hline
$\frac{1}{4}$&             5.10659e-01&             -&   7.83598e-02&               -&               4.17915e-02&   -\\
$\frac{1}{8}$&             2.60370e-01&       0.9718 &   2.08799e-02&               1.9080&          1.18158e-02 &       1.8225\\
$\frac{1}{16}$&            1.30809e-01&      0.9931&     5.30840e-03&               1.9758&          3.04107e-03&        1.9581\\
$\frac{1}{32}$&            6.54822e-02&       0.9983&    1.33278e-03&               1.9938&          7.65687e-04&         1.9898\\
\hline
  \end{tabular}
  \end{center}
\end{table}

\begin{table}[htpb]
         \centering
	\caption{ Convergence order of spatial error under PA5 parameter}
          \label{table15}
            \begin{center}
  \begin{tabular}{ccccccccc}
   \hline
 $h$&   $|\!|\!|{\mathbf u}(t_n)-{\mathbf u}_h^n|\!|\!|_V$&    R&   $||p(t_n)-p_h^n||$&         R&           $||T(t_n)-T_h^n||$&          R \\    \hline
$\frac{1}{4}$&             5.10316e-01&             -&   7.85361e-02&               -&               4.18521e-02&   -\\
$\frac{1}{8}$&             2.60308e-01&       0.9712 &   2.09566e-02&               1.9060&          1.18462e-02 &       1.8209\\
$\frac{1}{16}$&            1.30801e-01&      0.9928&     5.33026e-03&               1.9751&          3.05000e-03&        1.9575\\
$\frac{1}{32}$&            6.54811e-02&       0.9982&    1.33844e-03&               1.9937&          7.68020e-04&         1.9896\\
\hline
  \end{tabular}
  \end{center}
\end{table}

\begin{figure}[!htb]
\centering
\includegraphics[width=3.5in]{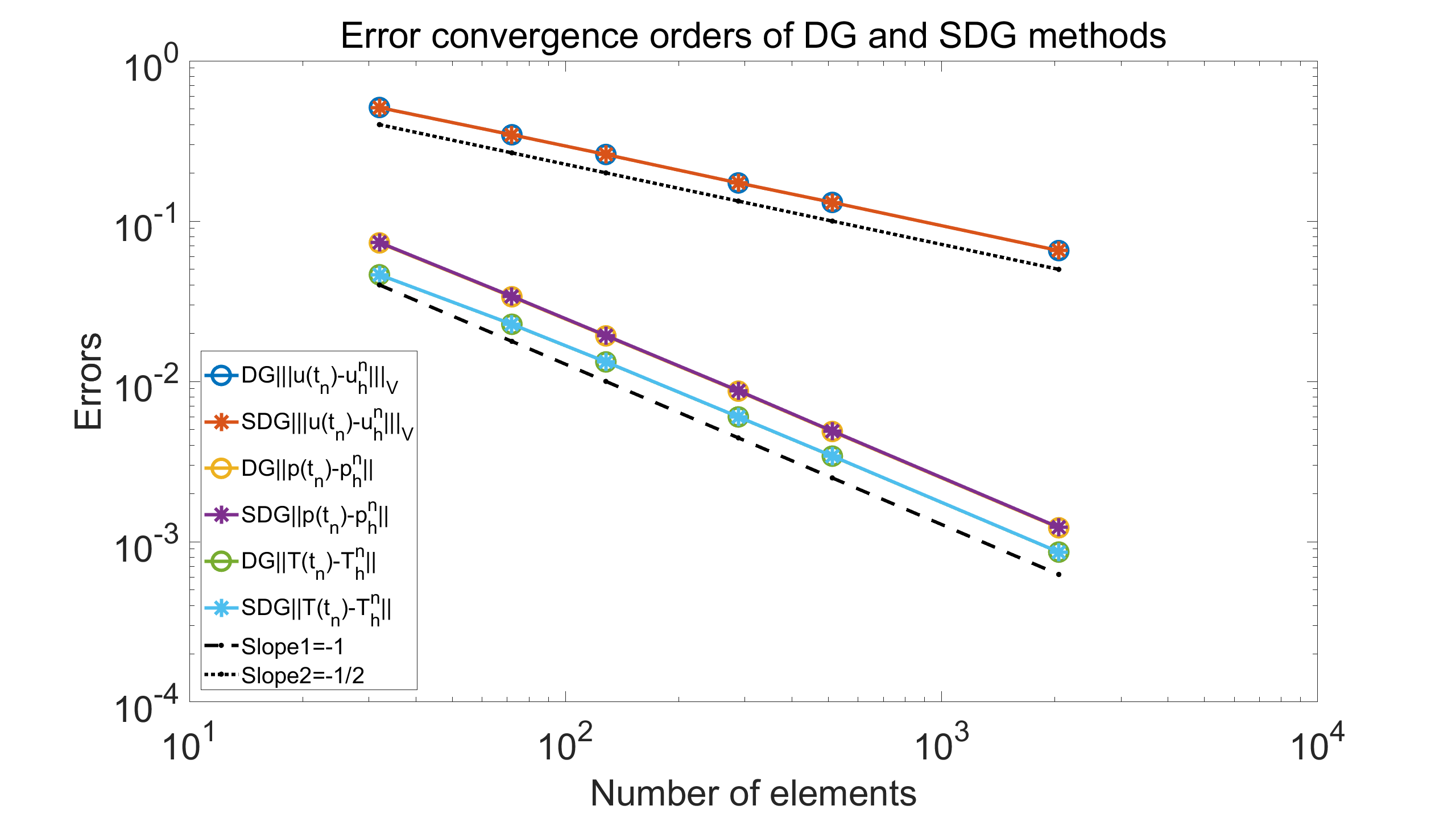}
\caption{\small{Error convergence orders of fully implicit DG and SDG methods}}\label{fig3-1}
\end{figure}

We also apply the fully implicit nonlinear numerical scheme (DG  method for spatial discretization and backward Euler method for temporal discretization) proposed in \cite{chen} for the system \eqref{e1}-\eqref{e3} in order to compare with
our method. The errors of our SDG method and fully implicit DG method are shown in Fig. \ref{fig3-1}, which indicates that both of these two methods have the similar error precision and convergence orders. In addition, we record the CPU running time of SDG method and fully implicit DG method with parameter choice PA1 in Table \ref{table17},
which shows that our method is more efficient than the fully implicit DG method.

\begin{table}[htpb]
\centering
\caption{ Comparison of computational efficiency for parameter choice PA1}
\label{table17}
\begin{center}
\begin{tabular}{ccccccccc}
\hline
           $h$&       SDG method CPU time (s)& DG method CPU time (s)     \\   \hline
$\frac{1}{8}$&  5.366 &21.505\\

$\frac{1}{16}$&            94.544      &503.269\\

$\frac{1}{32}$&           1515.250 &5834.460\\
 \hline
  \end{tabular}
  \end{center}
\end{table}

Finally, we test the other selection strategy of initial values, \textbf{Option 1}. Choose the initial time step size as $\tau_0=10^{-6}$ and the final time $t_f=1+10^{-6}$. Selection of other parameters is provided as mentioned before, the numerical test of spatial convergence orders for error. Table \ref{table18} shows that under the choice of parameter PA1, the errors of our SDG method can also reach the optimal convergence order.
\begin{table}[htpb]
         \centering
	\caption{Error convergence orders for parameter choice PA1 ($\tau_0=10^{-6}$)}
          \label{table18}
            \begin{center}
  \begin{tabular}{ccccccccc}
   \hline

  $h$&   $|\!|\!|{\mathbf u}(t_n)-{\mathbf u}_h^n|\!|\!|_V$&    R&   $||p(t_n)-p_h^n||$&         R&    $||T(t_n)-T_h^n||$&          R \\    \hline
$\frac{1}{4}$&             5.10973e-01&             -&         7.35407e-02&               -&              4.63424e-02&   -\\
$\frac{1}{8}$&             2.60436e-01&       0.9723 &         1.93515e-02&               1.9261&         1.32556e-02 &       1.8057\\
$\frac{1}{16}$&            1.30818e-01&             0.9934&    4.90206e-03&               1.9810&          3.42540e-03&                  1.9523\\
$\frac{1}{32}$&            6.54833e-02&             0.9984&    1.22958e-03&               1.9952&          8.63389e-04&                  1.9882\\
\hline
  \end{tabular}
  \end{center}
\end{table}

\section*{Acknowledgments}
This work was partially supported by the National Natural Science Foundation of China (No. 12371386, 12301465), the Research Foundation for Beijing University of Technology New Faculty (No. 006000514122516) and the Natural Science Foundation of Shandong Province (No. ZR2024MA056, ZR2022MA081).

\section*{Declarations}

\subsection*{Conflict of interest}
All authors declare that they have no conflict of interest.

\subsection*{Data availability}
All data generated or analysed during the current study are available from the corresponding author on reasonable request.


\end{document}